\pdfoutput=1
\documentclass[toc]{mathprint}
\usepackage{rutarmacros}
\usepackage{macros}
\usepackage{subcaption}
\title[Fibre stability of self-affine sets]{Fibre stability for dominated self-affine sets}

\author
  {Roope Anttila}
  {Research Unit of Mathematical Sciences, P.O.\ Box 8000, FI-90014 University of Oulu, Finland}
  {roope.anttila@oulu.fi}

\author
  {Alex Rutar}
  {Department of Mathematics and Statistics, University of Jyväskylä, P.O.\ Box 35 (MaD), FI-40014 University of Jyväskylä, Finland}
  {alex@rutar.org}

\begin{document}
\begin{abstract}
    Let $K$ be a planar self-affine set.
    Assuming a weak domination condition on the matrix parts, we prove for all backward Furstenberg directions $V$ that
    \begin{equation*}
        \max_{E\in\Tan(K)} \max_{x\in \pi_{V^\bot}(E)} \dimH (\pi_{V^\bot}^{-1}(x)\cap E) = \dimA K - \dimA \pi_{V^\bot}(K).
    \end{equation*}
    Here, $\Tan(K)$ denotes the space of weak tangents of $K$.
    Unlike previous work on this topic, we require no separation or irreducibility assumptions.
    However, if in addition the strong separation condition holds, then there exists a $V\in X_F$ so that
    \begin{equation*}
        \max_{x\in \pi_{V^\bot}(K)} \dimH (\pi_{V^\bot}^{-1}(x)\cap K) = \dimA K - \dimA \pi_{V^\bot}(K).
    \end{equation*}
    Our key innovation is an amplification result for slices of weak tangents via pigeonholing arguments.
\end{abstract}

\section{Introduction}
Since the seminal work of Marstrand \cite{zbl:0056.05504}, geometric properties of projections and slices of sets has been a fundamental and highly active research topic in fractal geometry.
Marstrand gave upper bounds for typical slices (in the sense of the Lebesgue measure on the space of directions) of planar Borel sets and showed that the Hausdorff dimensions of typical projections are as large possible.
These results have since been extended to higher dimensions and stronger results have been established under various structural assumptions, see for example \cite{arxiv:2308.08819,zbl:1414.28014,zbl:1426.11079,zbl:1430.11106,zbl:0819.28004,zbmath:7846313,zbmath:7810678}.

One caveat in Marstrand's theorems is that the sets of directions where the results hold are not explicit and bounding the dimensions of concrete slices and projections has remained a difficult problem.
In this context, the study of slices of planar self-affine sets has received much attention recently; for instance, we refer the reader to \cite{zbl:1371.28016, zbl:1455.28005,arxiv:2107.00983,doi:10.1017/etds.2023.117,zbl:07808129,zbl:1278.37032}.
It turns out that there is a natural connection between the \emph{Assouad dimension}, see \cref{ss:large-tangents}, of self-affine sets, and the dimensions of their slices and projections in certain directions called the (backward) Furstenberg directions.
These directions, which we denote by $X_F$, are essentially the ones in which deep iterates of the maps in the IFS are contracting as much as possible, see \cref{s:weak-domination} for a precise definition.
In general, it is expected that for self-affine sets $K$ there is a direction $V\in X_F$, such that
\begin{equation}\label{e:assouad-slice}
    \dimA K=\dimA \pi_{V^{\bot}}(K)+\max_{x\in\pi_{V^{\bot}}(K)}\dimA(\pi_{V^{\bot}}^{-1}(x)\cap K),
\end{equation}
where $\dimA K$ denotes the Assouad dimension of $K$ and $\pi_{V^{\bot}}$ is the orthogonal projection along the direction $V$.
This question was originally motivated in the case of diagonal matrices by Mackay \cite{zbl:1278.37032}, and for more general self-affine sets by the work of Bárány--Käenmäki--Rossi \cite{zbl:1455.28005}.
This also appears as a question explicitly in \cite[Question~17.5.1]{zbl:1467.28001}.
Equation \cref{e:assouad-slice} can be interpreted as a kind of stability under projections: the largest fibre stores the dimension lost in the projection.
With this motivation we call self-affine sets satisfying \cref{e:assouad-slice} \emph{fibre stable}.

Fibre stability has been established for various self-affine sets under a number of assumptions; for attractors of affine IFSs where each $A_i$ is a diagonal matrix see \cite{zbl:07808129,zbl:1371.28016,zbl:1278.37032,zbl:1305.28021}, and for more general self-affine sets see \cite{zbl:1455.28005,arxiv:2107.00983,doi:10.1017/etds.2023.117}.
In all of these works (with the exception of \cite{zbl:1305.28021} concerning Barański carpets), domination of the matrix parts of the IFS, see \cref{ss:dominated}, has played a crucial role.
Domination imposes a structure on the action of the matrix parts of the IFS on real projective space $\RP^1$ which is useful when analysing the directions in $X_F$.
Moreover, in all previous results, the authors have assumed that either the projections have Hausdorff dimension 1 (as follows from an irreducibility assumption on the matrix parts if the set has Hausdorff dimension at least 1) \cite{zbl:1455.28005,arxiv:2107.00983,doi:10.1017/etds.2023.117}, or are very well-behaved \cite{zbl:1371.28016,zbl:07808129,zbl:1278.37032,zbl:1305.28021}.
In particular, until the present work, \cite[Question~17.5.1]{zbl:1467.28001} was still open in its full generality for dominated self-affine sets satisfying the strong separation condition, even in the special case when all of the matrices are diagonal.

Our main result in this paper is a slicing theorem for dominated planar self-affine sets with no separation assumptions and no assumptions on the geometry of the projections.
In fact, for our purposes a slightly weaker variant of domination suffices: we say that an affine IFS is \emph{weakly dominated} if the matrix parts of the IFS can be decomposed into two subcollections, one of which is dominated and the other one of which consists of similarity maps which preserve the invariant directions, see \cref{d:weak-domination} for the precise definition.
Our proof technique is entirely self-contained and combines some basic geometric properties of matrix semigroups with pigeonholing arguments to obtain a type of amplification result for large slices.
Most notably, we do not depend on the deep work of \cite{zbl:1414.28014,zbl:1426.11079}, which has played a critical role in previous work.
Our results establish fibre stability for weakly dominated self-affine sets satisfying the strong separation condition in the plane and give the lower bound in \cref{e:assouad-slice} for the Assouad dimension with no separation assumptions.
We believe that the techniques in this paper will prove useful beyond the weakly dominated case.

\subsection{Large weak tangents and Assouad dimension}\label{ss:large-tangents}
Before stating our results let us set up some notation and recall some properties of the Assouad dimension.
The \emph{Assouad dimension} of a bounded set $F\subset \R^d$ is the number
\begin{align*}
    \dimA F=\inf\Big\{s\geq 0\colon &\exists C>0,\text{ such that }\forall 0<r\leq R\leq 1,\,x\in F\\
                                &N_r(F\cap B(x,R))\leq C\left(\frac{R}{r}\right)^s\Big\},
\end{align*}
where $N_r(E)$ denotes the smallest number of open balls of radius $r$ required to cover the set $E\subset \R^d$.
We note that at various points in the article we use alternate definitions, where $N_r(E)$ is replaced by the cardinality of the largest $r$-separated subset of $E$ or the dyadic covering number of $E$ at a scale $2^{-n}\approx r$.
Making these modifications does not affect the value of the Assouad dimension.

The Assouad dimensions of compact sets are closely related to the notion of a weak tangent.
For a closed set $F$, we denote by $\mathcal{K}(F)$ the set of all non-empty compact subsets of $F$ equipped with the \emph{Hausdorff metric} $d_{\mathcal{H}}$.
We then say that a set $E\in\mathcal{K}(B(0,1))$ is a \emph{microset} of $F$ if there exists a sequence $(x_n)_{n=1}^\infty\subset F$ and scales $(r_n)_{n=1}^\infty$ with $0<r_n\leq 1$ such that
\begin{equation*}
    E = \lim_{n\to\infty} r_n^{-1}(F - x_n) \cap B(0,1).
\end{equation*}
Moreover, we say that $E$ is a \emph{weak tangent} of $F$ if, in addition, $\lim_{n\to\infty} r_n=0$.
We denote the set of all microsets of $F$ by $\mathcal{G}_F$, and the set of all weak tangents of $F$ by $\Tan(F)$.
In general, $\Tan(F)\subset \mathcal{G}_F\subset\mathcal{K}(B(0,1))$.
We emphasize here that we do not permit rotations in our definition of a weak tangent.

A key observation is that the largest microset of a set exhibits substantially more regularity than the original set.
The proof of the following proposition is due to Furstenberg, but the explicit connection to Assouad dimension was first made in \cite{doi:10.1093/imrn/rnw336} and the amplification to Hausdorff content was noted in \cite{arxiv:2309.11971}.
\begin{proposition}[\cite{zbl:1154.37322}]\label{p:large-weak-tan}
    Let $F\subset\R^d$ be a compact set with $\eta = \dimA F$.
    Then $\dimH A\leq \eta$ for all $E \in\mathcal{G}_F$.
    Moreover, there is an $E\in\Tan(F)$ such that $\mathcal{H}^{\eta}_\infty(E) \geq 1$.
\end{proposition}
The main point is that a maximal microset of $F$ has $\dimH E = \dimA E = \dimA F$, which is a substantially gain in regularity over the original set.

\subsection{Main results}
Recall that an \emph{affine iterated function system (IFS)} is a finite collection $(T_i)_{i\in\mathcal{I}}$ of contracting invertible affine maps on $\R^2$, that is $T_i(x)=A_ix+b_i$, for each $i\in\mathcal{I}$, where $A_i$ is an invertible $2\times 2$ matrix with $\|A_i\|<1$ and $b_i\in\R^2$.
The \emph{attractor} of the IFS, which is also called the \emph{self-affine set}, is the unique non-empty and compact set $K$, which satisfies
\begin{equation*}
    K=\bigcup_{i\in\mathcal{I}}T_i(K).
\end{equation*}
We say that the IFS $(T_i)_{i\in\mathcal{I}}$ satisfies the \emph{strong separation condition}, if $T_i(K)\cap T_j(K)=\emptyset$ for all $i\ne j$.

Let $\RP^1$ denote the \emph{real projective space} of one dimensional subspaces of $\R^2$.
For $V,W\in\RP^1$ with $V\ne W$, we denote by $\pi_V^W\colon \R^2\to V$ the \emph{projection onto $V$ along $W$}, which is the unique linear map satisfying $\im(\pi_V^W)=V$ and $\ker(\pi_V^W)=W$ such that $\pi_V^W(v)=v$ for all $v\in V$.
If $W$ is not specified, then $\pi_V\colon \R^2\to V$ denotes the orthogonal projection onto $V$.

Our first result is a general slicing theorem for weak tangents, which surprisingly seems to have not been noticed before (in the context of self-affine sets) even though it is a straightforward consequence of Furstenberg's dimension conservation result \cite[Theorem 6.1]{zbl:1154.37322}.
The result also has a short elementary proof which we give in \cref{s:furstenberg} to keep the paper self-contained.
\begin{iproposition}\label{ip:large-tan-slice}
    Let $F\subset\R^2$ be an arbitrary non-empty compact set, and let $W\in \RP^1$ be arbitrary.
    Then there exists an $E\in\Tan(F)$ and $x\in\pi_W(E)$ such that
    \begin{equation*}
        \dimH(\pi_W^{-1}(x) \cap E) \geq \max\{\dimA F - \dimA \pi_W(F), 0\}.
    \end{equation*}
\end{iproposition}
The maximum is relevant since the Assouad dimension can in fact \emph{increase} under projection, even for self-similar sets: see, for instance, \cite[§3.1]{zbl:1305.28021}.
In general, it can happen for all directions $W$ that there exists a weak tangent $E$ such that $\dimA \pi_W(F) + \dimH(\pi_W^{-1}(x)\cap E)$ substantially exceeds the Assouad dimension of $F$.
As a simple example consider the following.
For each $W\in \RP^1$, let $C_W = (W\cup W^\perp)\cap B(0,1)$ denote the ``plus''-shaped set containing the origin oriented in direction $W$.
Let $(W_n)_{n=1}^\infty\subset\RP^1$ be dense, and define the set
\begin{equation*}
    C = \bigcup_{n=1}^\infty 2^{-(n+1)}\cdot C_{W_n}+(2^{-n},0).
\end{equation*}
Then $\dimA C = \dimH C = 1$, but for all $W\in\RP^1$, there exists $E\in\Tan(F)$, such that $\dimA \pi_W(C) + \dimH \pi_W^{-1}(0)\cap E = 2$.
Moreover, one cannot hope to improve \cref{ip:large-tan-slice} to slices of the original set; easy counterexamples are already given by function graphs of dimension strictly larger than $1$.

On the other hand, for self-affine sets we can say a lot more.
The following is the main result of this paper.
\begin{itheorem}\label{it:tan-slice-form}
    Let $(T_i)_{i\in\mathcal{I}}$ be a weakly dominated self-affine IFS with attractor $K$.
    Then the map $V\mapsto \dimA \pi_{V^{\bot}}(K)$ takes constant value $\eta\leq \dimA K$ on $X_F$.
    Moreover, for all $V\in X_F$,
    \begin{align*}
        \dimA K-\eta &= \max_{E\in\Tan(K)}\max_{x\in\pi_{V^{\bot}}(E)}\dimH(\pi_{V^{\bot}}^{-1}(x)\cap E)\\
                       &\geq \max_{x\in\pi_{V^\bot}(K)}\dimA(\pi_{V^\bot}^{-1}(x)\cap K)
    \end{align*}
    If in addition $(T_i)_{i\in\mathcal{I}}$ satisfies the strong separation condition, then
    \begin{equation*}
        \dimA K - \eta = \max_{V\in X_F}\max_{x\in\pi_{V^\bot}(K)}\dimH(\pi_{V^\bot}(x)^{-1}\cap K).
    \end{equation*}
\end{itheorem}
The proof of this result is split into multiple parts: \cref{t:weakly-dominated-assouad,p:sup-attained,c:real-slice-ub,c:real-slice}.
Also, for a reader only interested in the special case when the matrices are all diagonal, we give a condensed proof in \cref{as:diag}.

Let us make a few comments on \cref{it:tan-slice-form}.
\begin{enumerate}[nl]
    \item We require no assumptions concerning the projections of the self-affine set or irreducibility of the matrix parts.
        Most notably, our results also hold (and are new in this generality) for reducible self-affine sets, such as self-affine carpets satisfying the weak domination hypothesis.
    \item The conclusion concerning slices of weak tangents holds \emph{uniformly over all directions}: rather than the maximal value being attained at some direction in $X_F$, the maximum is attained in all directions in $X_F$ simultaneously.
        For slices of the set itself, it seems that there is no reason for this to be the case; see the proof of \cref{c:real-slice}.
        However, since there are not many tools to give non-trivial upper bounds for dimensions of all slices of a self-affine set in a given direction, coming up with counterexamples seems to be difficult.
    \item The results concerning slices of weak tangents and upper bounds on slices of $K$ hold with no separation assumptions at all; the planar separation is only required to ``pull back'' slices of the weak tangent to the original set (see \cref{c:real-slice} for the short proof).
    \item A somewhat weaker variant of the strong separation condition called the \emph{weak bounded neighbourhood condition} suffices; see \cref{d:wbnc}.
\end{enumerate}
In particular, this result substantially generalizes all previously known results concerning slices in backwards Furstenberg directions \cite{zbl:1371.28016, zbl:1455.28005,arxiv:2107.00983,doi:10.1017/etds.2023.117,zbl:07808129,zbl:1278.37032} and establishes fibre stability for weakly dominated and strongly separated self-affine sets.

Let us also emphasize that the value $\eta$ is the constant value of the Assouad dimension, rather than the Hausdorff dimension, of the relevant projections.
In \cite{zbl:1455.28005,arxiv:2107.00983,doi:10.1017/etds.2023.117}, this formula was written with Hausdorff dimension, but in those cases the assumptions implied that $\dimH \pi_{V^\bot}(K)=1$ for all $V\in X_F$.
Indeed, the following example follows from \cite[Theorem~2.13]{zbl:07808129}, using a similar construction as used in \cite[§2.5]{zbl:07808129} except with non-trivial fibres to guarantee that there is a symbolic slice with Assouad dimension 1.
\begin{proposition}[\cite{zbl:07808129}]\label{p:except}
    Let $\varepsilon>0$ be arbitrary.
    Then there is a planar dominated self-affine set $K$ with $\dimA K=2$ such that for all $V\in X_F$ and $x\in \pi_{V^\bot}(K)$, $\dimuB \pi_{V^\bot}(K)\leq\varepsilon$ and $\dimuB(\pi_{V^\bot}^{-1}(x)\cap K) \leq\varepsilon$.
\end{proposition}

This example also shows that in general, the separation assumption is needed for the second part of \cref{it:tan-slice-form} to hold as stated.
However, we are unsure whether or not the separation assumption is needed for the result to hold for the Assouad dimension of slices instead of Hausdorff dimension.
More precisely, we ask the following question.
\begin{question}\label{q:large-slice-ex}
    Let $(T_i)_{i\in\mathcal{I}}$ be a weakly dominated self-affine IFS with attractor $K$ and let $\eta$ be the constant value of the map $V\mapsto \dimA \pi_{V^{\bot}}(K)$ on $X_F$.
    Is it true that
    \begin{equation*}
        \dimA K - \eta = \max_{V\in X_F}\max_{x\in\pi_{V^\bot}(K)}\dimA(\pi_{V^\bot}^{-1}(x)\cap K)?
    \end{equation*}
\end{question}

To conclude the introduction, let us note two direct applications of \cref{it:tan-slice-form}.
First, because of the regularity of the space of weak tangents, assuming separation \cref{it:tan-slice-form} is equivalently a result concerning the dimensions of \emph{tubes}; see \cref{c:tube-bound}.
In particular, we are able to complete a partial result due to Fraser \& Jordan concerning certain self-affine carpets with no grid structure.
Let $0<\alpha<\beta<1$ and consider the self-affine system defined by maps $T_i(x,y)=(\beta x, \alpha y) + (b_i, a_i)$ for $0\leq b_i \leq 1-\beta$ and $0\leq a_i \leq 1-\alpha$, and let $\nu$ denote the uniform self-similar measure associated with the projected IFS defined by maps $(\beta x+b_i)_{i\in\mathcal{I}}$.
Finally, let $s=\dim_\infty\nu$ denote the Frostman dimension of $\nu$.
The following result generalizes \cite[Theorem~2.2]{zbl:1371.28016}; for the proof, along with more careful exposition, see \cref{ss:geom-bound}.
\begin{icorollary}\label{ic:tubes}
    Let $K$ be the self-affine set corresponding to the IFS $(T_i)_{i\in\mathcal{I}}$ defined above, with parameters $0<\alpha<\beta<1$.
    Assume moreover that the $T_i((0,1)^2)\cap T_j((0,1)^2)=\varnothing$ for all $i\neq j$.
    Let $s$ denote the Frostman dimension of $\nu$.
    Then
    \begin{equation*}
        \dimA K = \dimA \pi(K) + \frac{\log m\beta^s}{\log(1/\alpha)}.
    \end{equation*}
\end{icorollary}

Finally, we note an application to conformal Assouad dimension, which follows from \cref{it:tan-slice-form} combined with a strong projection theorem for Assouad dimension due to Orponen \cite{zbl:1465.28008}.
We recall that the conformal Assouad dimension is defined by
\begin{equation*}
    \CdimA X\coloneqq \inf\left\{\dimA f(X)\colon f\text{ is a quasisymmetry}\right\}.
\end{equation*}
We refer the reader to \cref{ss:conformal} for more background.
Also, recall that a self-affine set is \emph{irreducible} if no linear subspace in $\RP^1$ is preserved by all of the linear parts of the affine maps in the IFS.
Unlike similar results which have previously appeared (such as \cite[Theorem~3.2]{zbl:1455.28005} and \cite[Theorem~B]{doi:10.1093/imrn/rnw336}), we require no separation assumptions either in the plane or in the projection.
The proof can be found in \cref{ss:conformal}.
\begin{icorollary}\label{ic:conformal}
    Let $K$ be a weakly dominated and irreducible self-affine set.
    If $\dimA K<1$, then $\CdimA K=0$, and if $\dimA K\geq 1$, then $K$ is minimal for conformal Assouad dimension.
\end{icorollary}

\subsection{Outline of paper}
In \cref{s:furstenberg}, we establish some preliminaries concerning microsets and weak tangents, and in particular  in \cref{l:h-micro-hausdorff} we give a short proof of the discretized variant of Furstenberg's microset existence argument.
We also give the self-contained proof of \cref{ip:large-tan-slice}.
Next, in \cref{s:weak-domination} we establish some preliminaries concerning weak domination; the results stated here are relatively standard and are mostly drawn either from \cite{zbl:1452.37039} or recent papers concerning self-affine sets.

The heart of the paper is \cref{s:proof}, where we establish the main slicing result for weak tangents, stated in \cref{t:weakly-dominated-assouad}.
The key innovation is a combination of \cref{l:h-micro-hausdorff} to show the existence of microsets with large covering numbers across arbitrary sequences of scales with a delicate pigeonholing argument to find a collection of well-aligned copies of approximations of maximal weak tangents of projections inside the self-affine set.
Using the self-affine structure, this configuration, which is a subset of a thin tube in the direction of the slice, can be pushed to a product-like structure inside some well chosen cylinder.
This approach is made formal in \cref{t:product-tangent}; a more precise (but still informal) discussion of the proof can also be found immediately preceding \cref{t:product-tangent}.

We note for the reader only interested in the special case of \cref{it:tan-slice-form} for diagonal matrices, since the geometry of the matrix semigroup is very simple in this case, \cref{s:weak-domination} and \cref{ss:geom} can be skipped entirely and a condensed proof can be found in \cref{as:diag}.

Finally, in \cref{s:consequences}, we complete the remaining minor components of the proof of \cref{it:tan-slice-form}.
We then discuss the application to certain diagonal self-affine sets and prove \cref{ic:tubes}, and the application to the conformal Assouad dimension of weakly dominated and irreducible self-affine sets with no separation assumptions, and prove \cref{ic:conformal}.

\subsection{Notation}
We use $\langle \cdot, \cdot\rangle$ to denote the standard dot product on $\R^2$.

For a linear subspace $V$ of $\R^2$ and a linear map $T\colon V\to \R^2$, we denote by $\|T\|$ the \emph{operator norm} of $T$, that is
\begin{equation*}
    \|T\|=\max_{v\in V\setminus\{0\}}\frac{\|Tv\|}{\|v\|},
\end{equation*}
where $\|\cdot\|$ is the standard Euclidean norm on $\R^2$
Given a subspace $V\subset \R^2$, and a $2\times 2$ matrix $A$, which we interpret as a linear map from $\R^2$ to itself, we denote the \emph{restriction of $A$ to $V$} by $A|V\colon V\to\R^2$.

We will sometimes make use of the following asymptotic notation.
Given a set $A$ and functions $f,g\colon A\to \R$ we write $f\gtrsim g$, if there is a constant $C$, such that $f(a)\geq Cg(a)$, for all $a\in A$.

\section{Amplifying dimension and slicing weak tangents}\label{s:furstenberg}
In this section, we introduce the techniques we use to bound the Assouad dimension from below, the most important of which is a discretized variant of Furstenberg's well known construction for measures to show the existence of microsets with uniformly large branching over arbitrarily long sequences of scales.
This construction plays a crucial role in the proofs of our main results, and enables us to give a short and self-contained proof of \cref{ip:large-tan-slice}.

Let us start by introducing some basic notation concerning dyadic cubes.
Fix $d\in\N$.
Let $\mathcal{D}=\bigcup_{n=0}^\infty\mathcal{D}_n$ denote the set of closed dyadic cubes, where $\mathcal{D}_n$ denotes the subset of dyadic cubes with side-length $2^{-n}$.
Given $Q\in\mathcal{D}$, let: $\psi_Q\colon Q\to[0,1]^d$ denote the unique surjective homothety mapping $Q$ to $[0,1]^d$.
For $n\in\N$, and a bounded set $K\subset \R^d$ we let $N_n$ denote the level $n$ dyadic covering number of $K$, that is
\begin{equation*}
    N_n(K)=\#\left\{Q\in\mathcal{D}_n:Q\cap K\neq\varnothing\right\}.
\end{equation*}
We note that the covering numbers have the property that for any $Q\in\mathcal{D}_m$ and $n\in\N$,
\begin{equation*}
    N_{m+n}(K\cap Q)=N_n(\psi_{Q}(K\cap Q))=N_n(\psi_{Q}(K)\cap Q_0).
\end{equation*}
This simple property will be used throughout the rest of the paper without further reference.

\subsection{Coarse microsets}
Let $K_1,K_2\subset\R^d$ be a non-empty compact set.
Let $p_{\mathcal{H}}$ denote the one-sided Hausdorff metric
\begin{equation*}
    p_{\mathcal{H}}(K_1;K_2)=\inf\{\delta\geq 0:K_1\subset K_2^{(\delta)}\}.
\end{equation*}
Here, $K_2^{(\delta)}$ denotes the open $\delta$-neighbourhood of $K_2$.

Now let $K$ be a non-empty compact set.
We say that a non-empty compact set $E$ is a \emph{coarse microset} of $K$ if there is a sequence of expansion ratios $\lambda_n\geq 1$, points $x_n\in K$, and a bi-Lipschitz map $f\colon\R^d\to\R^d$ such that
\begin{equation*}
    \lim_{n\to\infty} p_{\mathcal{H}}\bigl(f(E);\lambda_n (K-x_n)\bigr)=0.
\end{equation*}
The following lemma is standard.
\begin{lemma}\label{l:coarse-micro}
    Let $K\subset\R^d$ be non-empty and compact.
    Then $\dimA K\geq\dimA E$ for any coarse microset $E$ of $K$.
\end{lemma}
In order to lower bound the Assouad dimension of a coarse microset, we also note the following standard lemma which follows by semi-continuity of dyadic covering numbers.
\begin{lemma}\label{l:lower-box}
    Let $F\subset\R^d$ be non-empty and compact and let $(F_n)_{n=1}^\infty$ be a sequence of non-empty compact sets such that
    \begin{equation*}
        \lim_{n\to\infty} d_{\mathcal{H}}(F,F_n)=0.
    \end{equation*}
    Suppose moreover that there is an unbounded sequence of natural numbers $(m_n)_{n=1}^{\infty}$, such that
    \begin{equation*}
        N_{k}(F_n)\gtrsim 2^{ks}
    \end{equation*}
    for all $0\leq k \leq m_n$.
    Then $\dimlB F\geq s$.
\end{lemma}

\subsection{Dyadic cubes and weak tangents}\label{ss:furstenbergs-lemma}
We now demonstrate the existence of minisets with uniformly large branching over arbitrarily large sequences of levels.
This is Furstenberg's well-known pigeonholing construction for measures; see, for instance, \cite[Lemma~2.4.4]{zbl:1390.28012} or \cite[Theorem~5.1.3]{zbl:1467.28001}.
Note that \cref{e:branching} is a Frostman-type condition for the measure on $Q$ which is uniformly distributed on $Q\cap K$ at level $m$.
\begin{lemma}\label{l:h-micro-hausdorff}
    Let $K\subset[0,1]^d$ be a non-empty compact set.
    Let $0<s<t$, $\ell\in\N$, and $k\in\N$ with $k\geq\ell$.
    Suppose there is $m\geq\frac{k d}{t-s}$ so that
    \begin{equation}\label{e:top-count}
        N_{m}(K)\geq 2^{m t}.
    \end{equation}
    Then there is a $0\leq p\leq m-k$ and a dyadic cube $Q\in\mathcal{D}_p$ so that for all $0\leq n\leq \ell$ and $Q\supset Q'\in\mathcal{D}_{p+n}$,
    \begin{equation}\label{e:branching}
        \frac{N_m(K\cap Q')}{N_m(K\cap Q)} \leq 2^{-ns}.
    \end{equation}
\end{lemma}
\begin{proof}
    If $Q_0=[0,1]^d$ satisfies the branching condition \cref{e:branching}, we are done.
    Otherwise, there is $1\leq\ell_1\leq \ell$ and a $Q_1\in\mathcal{D}_{\ell_1}$ so that
    \begin{equation}\label{e:sub-count}
        N_{m}(K\cap Q_1)>2^{m t}2^{-\ell_1 s}.
    \end{equation}
    Repeating the above argument for each $j\geq 1$ with $K\cap Q_j$ in place of $K$ and $m-\ell_1-\cdots-\ell_j$ in place of $m$, either there is some count $q$ with $\ell_1+\cdots+\ell_q\leq m-k$ such that the dyadic cube $Q_q\in\mathcal{D}_{\ell_1+\cdots+\ell_q}$ satisfies the branching condition, or $m\geq\ell_1+\cdots+\ell_q>m-k$.
    Suppose for contradiction that the latter situation occurs.
    Then
    \begin{equation*}
        2^{k d}\geq N_{m}(K\cap Q_k)>2^{m t}2^{-(\ell_1+\cdots+\ell_k) s}\geq 2^{m(t-s)}.
    \end{equation*}
    Rearranging, $k d > m (t-s)$ which contradicts the choice of $m$.
\end{proof}
By pigeonholing, we recover the following slightly weaker version of the conclusion which we find somewhat more convenient to use.
\begin{corollary}\label{l:h-micro}
    Let $K\subset[0,1]^d$ be a non-empty compact set.
    Let $0<s<t$, $\ell\in\N$, and $k\in\N$ with $k\geq\ell$.
    Suppose there is $m\geq\frac{k d}{t-s}$ so that
    \begin{equation*}
        N_{m}(K)\geq 2^{m t}.
    \end{equation*}
    Then there is a $0\leq p\leq m-k$ and a dyadic cube $Q\in\mathcal{D}_p$ so that
    \begin{equation*}
        N_{p+n}(K\cap Q)\geq 2^{ns}
    \end{equation*}
    for all $0\leq n\leq \ell$.
\end{corollary}
Combining \cref{l:h-micro} with the definition of the Assouad dimension yields the following.
\begin{corollary}\label{c:weak-tan}
    Let $K\subset[0,1]^d$ be a non-empty compact set with $\dimA K=\eta$.
    Then there is a sequence of dyadic cubes $(Q_m)_{m=1}^\infty\subset\mathcal{T}$ with $\diam Q_m$ decreasing to $0$ so that
    \begin{equation*}
        N_{n}(\psi_{Q_m}(K)\cap Q_0)\geq 2^{n\left(\eta-\frac{1}{m}\right)}
    \end{equation*}
    for all $0\leq n\leq m$.
\end{corollary}

These simple lemmas enable us to give a short and elementary proof of \cref{ip:large-tan-slice} which we restate here for convenience.
\begin{restatement}{ip:large-tan-slice}
    Let $F\subset\R^2$ be an arbitrary non-empty compact set, and let $W\in \RP^1$ be arbitrary.
    Then there exists an $E\in\Tan(F)$ and $x\in\pi_W(E)$ such that
    \begin{equation*}
        \dimH(\pi_W^{-1}(x) \cap E) \geq \max\{\dimA F - \dimA \pi_W(F),0\}.
    \end{equation*}
\end{restatement}
\begin{proof}
    The result is clearly true of $\dimA\pi_W(F) > \dimA F$, so we may assume otherwise.
    Denote by $s=\dimA F$ and $\eta=\dimA \pi_W(F)$, and let $Q_0=[0,1]^2$.
    By rotating the set $F$ if necessary, we will assume that $\pi_W=\pi$ is the projection on the $x$-axis.
    Let $n\in\N$ and take $m\geq2n^2$ and by \cref{c:weak-tan} pick $Q\in \mathcal{D}$, such that
    \begin{equation*}
        N_{m}(\psi_Q(F)\cap Q_0)\geq 2^{m(s-\frac{1}{n})}.
    \end{equation*}
    Let $\mathcal{P}_m$ denote the partition of the unit square into congruent tubes of width $2^{-m}$ and height one.
    Note that by the definition of the Assouad dimension, the set $\pi(\psi_{Q}(F)\cap Q_0)$ intersects at most $2^{m(\eta+\frac{1}{n})}$ dyadic intervals of length $2^{-m}$ and therefore, $\psi_Q(F\cap Q)$ intersects at most $2^{m(\eta+\frac{1}{n})}$ tubes of width $2^{-m}$.
    By the pigeonhole principle, there is $P_n\in\mathcal{P}_m$, such that
    \begin{equation*}
        N_{m}(\psi_Q(F)\cap P_n)\geq 2^{m(s-\eta-\frac{2}{n})}.
    \end{equation*}
    Now apply \cref{l:h-micro}, with $t=s-\frac{2}{n}$, $s=t-\frac{1}{n}$ and $k=\ell=n$ to find a dyadic cube $Q_n\in\mathcal{D}_p$ for $0\leq p\leq m-n$, such that
    \begin{equation}\label{e:large-tube}
        N_{k}(\psi_{Q_n}(\psi_Q(F)\cap P_n)\cap Q_0)\geq 2^{k(s-\eta-\frac{3}{n})},
    \end{equation}
    for all $0\leq k \leq n$.
    Note that $\psi_{Q_n}(P_n)$ is a tube of width $2^{-(m-p)}\leq 2^{-n}$, so by passing to a subsequence, there exist compact sets $E\in\Tan(F)$ and $A\subset E$ and a point $x\in\pi(E)$, such that $\psi_{Q_n}(\psi_{Q}(F))\cap Q_0\to E$ and $\psi_{Q_n}(\psi_Q(F)\cap P_n)\cap Q_0\to A\subset \pi^{-1}(x)\cap E$.
    Therefore, by \cref{e:large-tube} and \cref{l:lower-box},
    \begin{equation*}
        \dimlB(\pi^{-1}(x)\cap E)\geq s-\eta.
    \end{equation*}
    Using \cref{p:large-weak-tan} to pass again to a weak tangent of $\pi^{-1}(x)\cap E$ yields the desired result for Hausdorff dimension.
\end{proof}

\section{Weak domination in matrix semigroups}\label{s:weak-domination}
Let $\MM_2$ denote the space of $2\times 2$ real matrices and $\GL_2$ denote the group of invertible matrices in $\MM_2$.
Let $\mathcal{I}$ be a finite index set and let $\A=(A_i)_{i\in\mathcal{I}}$ be a tuple of matrices in $\GL_2$.
In the theory of self-affine sets, the action of the matrix semigroup generated by the tuple of the linear parts of the affine maps in the IFS, plays an important role.
In this section we describe this action in detail for weakly dominated matrices.

When studying matrix semigroups arising from affine IFSs, it is often useful to phrase the results with respect to the underlying symbolic space.
We call the symbol $\varnothing$ the \emph{empty word} and let $A_{\varnothing}=\mathrm{Id}$.
For $n\in\N$, let $\mathcal{I}^n$ denote the \emph{words of length $n$} generated by $\mathcal{I}$ and $\mathcal{I}^*=\bigcup_{n=0}^{\infty}\mathcal{I}^n$ denote the collection of all \emph{finite words}, where $\mathcal{I}^0=\{\varnothing\}$.
We call $\Sigma(\mathcal{I})=\mathcal{I}^{\N}$ the \emph{symbolic space} associated with $\mathcal{I}$.
If $\mathcal{I}$ is clear from the context, we may drop it from the notation and simply use the notation $\Sigma$ for $\Sigma(\mathcal{I})$.
We use the notation $\mtt{i}$ for words in both $\mathcal{I}^*$ and $\Sigma(\mathcal{I})$, that is $\mtt{i}=i_1i_2\cdots i_n$ and $\mtt{i}=i_1i_2\cdots$, respectively.
For $\mtt{i}\in\mathcal{I}^*$ we denote by $|\mtt{i}|$ the \emph{length} of $\mtt{i}$, which is the unique integer $n$, such that $\mtt{i}\in\mathcal{I}^n$.
For $\mtt{i}\in\Sigma(\mathcal{I})$ and $n\in\N$, we let $\mtt{i}|_n\coloneqq i_1i_2\cdots i_n$ denote the \emph{restriction of $\mtt{i}$} onto the first $n$ symbols.
For $\mtt{i}=i_1i_2\cdots i_n\in\mathcal{I}^*$, we let $\mtt{i}^-=i_1i_2\cdots i_{n-1}$.
Given any collection of functions $(F_i)_{i\in\mathcal{I}}$ or real numbers $(a_i)_{i\in\mathcal{I}}$, we denote for $\mtt{i}\in\mathcal{I}^n$,
\begin{align*}
    &F_{\mtt{i}} = F_{i_1}\circ F_{i_2}\circ\cdots\circ F_{i_n},\\
    &a_{\mtt{i}} = a_{i_1}a_{i_2}\cdots a_{i_n}.
\end{align*}

For a matrix tuple $\A$ we let
\begin{equation*}
    \mathfrak{R}(\A)=\{A\in\overline{\{cA_{\mtt{i}}\colon c\in\R\,\text{ and }\mtt{i}\in\mathcal{I}^*\}}\colon \mathrm{rank}(A)=1\}.
\end{equation*}
For us, there are two important sets of directions in $\RP^1$, namely the sets
\begin{align*}
    Y_F(\A)&\coloneqq\{\im(A)\in\RP^1\colon A\in\mathfrak{R}(\A)\},\\
    X_F(\A)&\coloneqq\{\im(A)\in\RP^1\colon A\in\mathfrak{R}(\A^{-1})\},
\end{align*}
where $\A^{-1}\coloneqq(A_i^{-1})_{i\in\mathcal{I}}$.
We call these sets the \emph{forward and backward Furstenberg directions}, respectively.
It is immediate from the definitions that $Y_F(\A)=X_F(\A^{-1})$.

\subsection{Dominated matrices}\label{ss:dominated}
Recall that $\A$ is dominated if there exist constants $0<\tau<1$ and $c>0$ such that
\begin{equation*}
    \alpha_2(A_{\mtt{i}})\leq c\tau^{|\mtt{i}|}\alpha_1(A_{\mtt{i}}),
\end{equation*}
for all $\mtt{i}\in\mathcal{I}^*$.
By \cite{zbl:1181.37032}, this is equivalent to the existence of a \emph{strongly invariant multicone} $\mathcal{C}\subset \RP^1$, meaning that $\mathcal{C}$ is a finite union of closed projective intervals satisfying $A_i(\mathcal{C})\subset \mathcal{C}^{\circ}$, for all $i\in\mathcal{I}$.

The Furstenberg directions for dominated tuples have a useful symbolic representation which we describe next.
For $A\in \MM_2$, we denote by $\alpha_1(A)\geq \alpha_2(A)$ the \emph{singular values} of $A$.
Formally, these are the square roots of the non-negative eigenvalues of the positive definite matrix $A^{\top}A$, and geometrically, they correspond to the lengths of the semiaxes of the ellipse $A(B(0,1))$.
The \emph{right singular vectors} of $A$ are eigenvectors $\eta_1(A)$ and $\eta_2(A)$ of $A^{\top}A$ corresponding to the eigenvalues $\alpha_1(A)$ and $\alpha_2(A)$, respectively.
If $\alpha_1(A)>\alpha_2(A)$, which is the case for all matrices in the semigroup generated by a dominated tuple, then these vectors are unique up to a change of sign.
For $\mtt{i}\in\mathcal{I}^n$, we write $\overleftarrow{\mtt{i}}=i_ni_{n-1}\ldots i_1$, and
\begin{align*}
    &A_{\overleftarrow{\mtt{i}}}^{-1}=(A_{\overleftarrow{\mtt{i}}})^{-1}=A_{i_1}^{-1}\cdots A_{i_n}^{-1}.
\end{align*}
We emphasize that $A_{\mtt{i}}^{-1}=A_{i_n}^{-1}\cdots A_{i_1}^{-1}$ denotes the inverse matrix of $A_{\mtt{i}}$, so in general, $A_{\mtt{i}}^{-1}\ne A_{\overleftarrow{\mtt{i}}}^{-1}$.
For $\mtt{i}\in\mathcal{I}^*$, we let
\begin{align*}
    \vartheta_1(\mtt{i})&=\langle A_{\mtt{i}}\eta_1(A_{\mtt{i}})\rangle\\
    \vartheta_2(\mtt{i})&=\langle A_{\overleftarrow{\mtt{i}}}^{-1}\eta_1(A_{\overleftarrow{\mtt{i}}}^{-1})\rangle,
\end{align*}
where $\langle v\rangle\in\RP^1$ denotes the line spanned by $v\in\R^2$.
The geometric interpretation therefore is that $\vartheta_1(\mtt{i})$ and $\vartheta_2(\mtt{i})$ are the lines spanned by the longer semiaxis of the ellipses $A_{\mtt{i}}(B(0,1))$ and $A_{\overline{\mtt{i}}}^{-1}(B(0,1))$, respectively.
For $\mtt{i}\in\Sigma$ and $k\in\{1,2\}$, we define
\begin{equation*}
    \overline{\vartheta_k}(\mtt{i})=\lim_{n\to\infty}\vartheta_k(\mtt{i}|_n),
\end{equation*}
whenever the limit exists.
It turns out that for dominated tuples, the limit always exists and the convergence is uniform, which means that one can think of $\overline{\vartheta_1}$ and $\overline{\vartheta_2}$ as projections from the symbolic space to $Y_F(\A)$ and $X_F(\A)$, respectively.
The proof of the following lemma can be found, for instance, in \cite[Lemmas~2.2 and 2.3]{doi:10.1017/etds.2023.117}.
\begin{lemma}\label{l:dominated}
    If $\A$ is dominated and $\mathcal{C}\subset \RP^1$ is a strongly invariant multicone for $\A$, then for $k\in\{1,2\}$:
    \begin{enumerate}[nl]
        \item the limit $\overline{\vartheta_k}(\mtt{i})=\lim_{n\to\infty}\vartheta_k(\mtt{i}|_n)$ exists for every $\mtt{i}\in\Sigma$ and the convergence is uniform;
        \item the map $\overline{\vartheta_k}(\mtt{i})\colon \Sigma\to \RP^1$ is Hölder continuous;
        \item\label{im:acc-proj} the set $\overline{\vartheta_k}(\Sigma)$ is compact and contains the accumulation points of $\{\overline{\vartheta_k}(\mtt{i})\colon \mtt{i}\in\mathcal{I}^*\}$;
        \item\label{im:shift-inv} $A_{\mtt{i}}\overline{\vartheta_1}(\mtt{j})=\overline{\vartheta_1}(\mtt{i}\mtt{j})$ and $A_{\overleftarrow{\mtt{i}}}^{-1}\overline{\vartheta_2}(\mtt{j})=\overline{\vartheta_2}(\mtt{i}\mtt{j})$, for all $\mtt{i}\in\mathcal{I}^*$ and $\mtt{j}\in\Sigma$;
        \item $Y_F(\A)=\overline{\vartheta_1}(\Sigma)\subset \mathcal{C}^{\circ}$ and $X_F(\A)=\overline{\vartheta_2}(\Sigma)\subset \RP^1\setminus\mathcal{C}$.
    \end{enumerate}
\end{lemma}
Another useful property of dominated tuples is that the singular values of the matrices in the semigroup are determined by restricting the matrices to suitable subspaces.
The next lemma follows from \cite[Lemma~2.8]{arxiv:2107.00983} by observing that $X_F(\A)^{\bot}=Y_F(\A^{\top})$.
\begin{lemma}\label{l:dominated-restrictions}
    If $\A$ is dominated, then there exists a constant $D\geq 1$ such that
    \begin{equation*}
        \|A_{\mtt{i}}|Y\|\leq \alpha_1(A_{\mtt{i}})\leq D\|A_{\mtt{i}}|Y\|,
    \end{equation*}
    for all $\mtt{i}\in\mathcal{I}^*$ and $Y\in Y_F(\A)$.
    Furthermore, if $V\in X_F(\A)$ and $\mtt{i}\in\Sigma$ is such that $V=\overline{\vartheta_2}(\mtt{i})$, then
    \begin{equation*}
        D^{-1}\|A_{\overleftarrow{\mtt{i}|_n}}|V\|\leq \alpha_2(A_{\overleftarrow{\mtt{i}|_n}})\leq \|A_{\overleftarrow{\mtt{i}|_n}}|V\|,
    \end{equation*}
    for all $n\in\N$.
\end{lemma}

\subsection{Weakly dominated matrices}
For our purposes, a slightly weaker variant of domination is sufficient.
We call a tuple $\A=(A_i)_{i\in\mathcal{I}}$ \emph{strongly conformal} if there exists a \emph{conjugation matrix} $M\in \GL_2$ such that for all $i\in\mathcal{I}$,
\begin{equation*}
    A_i=a_iM O_i M^{-1},
\end{equation*}
for some $0<a_i<1$ and $O_i\in \OO_2$, where $\OO_2$ denotes the subgroup of orthogonal matrices in $\GL_2$.
\begin{definition}\label{d:weak-domination}
    We say that $\A=(A_i)_{i\in\mathcal{I}}$ is \emph{weakly dominated} if it can be decomposed into two sets $\A_e$ and $\A_h$ such that $\A_e$ is strongly conformal and $\A_h$ is non-empty and has a strongly invariant multicone $\mathcal{C}$ such that $A\mathcal{C}=\mathcal{C}$ for all $A\in \A_e$.
\end{definition}
By a suitable change of coordinates determined by the conjugation matrix $M$, we may assume without loss of generality that each $A_i\in\A_e$ is of the form $a_iO_i$, for some $0<a_i<1$ and $O_i\in \OO_2$.
It was shown in \cite[Corollary 2.5]{zbl:1452.37039} that domination can be restated in terms \emph{almost multiplicativity} of the associated semigroup in the following sense.
\begin{lemma}\label{l:almost-multiplicative}
    A tuple $\A$ is either weakly dominated or strongly conformal if and only if there exists a constant $C>0$ such that
    \begin{equation*}
        C\|A_{\mtt{i}}\|\|A_{\mtt{j}}\|\leq \|A_{\mtt{i}\mtt{j}}\|\leq \|A_{\mtt{i}}\|\|A_{\mtt{j}}\|,
    \end{equation*}
    for all $\mtt{i},\mtt{j}\in\mathcal{I}^*$.
\end{lemma}

Next we will show that the Furstenberg directions of a weakly dominated tuple $\A$ are determined by a canonical dominated tuple $\OA$.
The ideas are essentially from \cite{zbl:1452.37039}, but we rewrite them with somewhat different notation more suitable for our purposes.
For the remainder of this section, fix a weakly dominated tuple $\A=(A_i)_{i\in\mathcal{I}}$ and let $\mathcal{I}_e=\{i\in\mathcal{I}\colon A_i\in \A_e\}$ and $\mathcal{I}_h=\{i\in\mathcal{I}\colon A_i\in \A_h\}$.
We define an equivalence relation on $\{\mtt{j}i\mtt{i}\in\mathcal{I}^*\colon i\in\mathcal{I}_h,\,\mtt{j},\mtt{i}\in\mathcal{I}_e^*\}$ by saying that $\mtt{j}_1i_1\mtt{i}_1\sim  \mtt{j}_2i_2\mtt{i}_2$ if and only if $O_{\mtt{j}_1}=O_{\mtt{j}_2}$ and $O_{\mtt{i}_1}=O_{\mtt{i}_2}$.
We denote the equivalence class of $\mtt{j}i\mtt{i}$ under this equivalence relation by $[\mtt{j}i\mtt{i}]$ and the collection of all equivalence classes by $\Lambda$.
It follows from \cite[Theorem 2.1 and Lemma 3.7]{zbl:1452.37039} that the sub-semigroup generated by the set $\{O_i\colon i\in\mathcal{I}_e\}$ is finite and therefore $\Lambda$ is a finite set.
Note that for any $\lambda=[\mtt{j}i\mtt{i}]\in\Lambda$ the matrix
\begin{equation*}
    \overline{A}_{\lambda}\coloneqq O_{\mtt{j}}A_{i}O_{\mtt{i}},
\end{equation*}
is well defined.

Every word $\mtt{i}\in\mathcal{I}^*\setminus\mathcal{I}_e^*$ can be uniquely decomposed as $\mtt{i}=\mtt{i}_0i_1\mtt{i}_2i_2\cdots \mtt{i}_ki_k$, for some $k\in\N$, where $\mtt{i}_j\in\mathcal{I}_e^*$ and $i_j\in \mathcal{I}_h$, for all $j=0,\ldots,k$.
Therefore, we may define a mapping $\mtt{i}\mapsto [\mtt{i}]$ from $\mathcal{I}^*\setminus \mathcal{I}_e^*$ to $\Lambda^*$ by setting
\begin{equation*}
    [\mtt{i}]=[\mtt{i}_0i_1\mtt{i}_1][\varnothing i_2\mtt{i}_2]\cdots [\varnothing\mtt{i}_ki_k].
\end{equation*}
Given $[\mtt{i}]\in\mathcal{I}^*\setminus \mathcal{I}_e^*$, we set
\begin{equation}\label{e:products}
    \begin{aligned}
        \overline{A}_{[\mtt{i}]}\coloneqq{}& \overline{A}_{[\mtt{i}_0i_1\mtt{i}_1]}\overline{A}_{[\varnothing i_2\mtt{i}_2]}\cdots \overline{A}_{[\varnothing i_k\mtt{i}_k]}\\
        ={}&O_{\mtt{i}_0}A_{i_1}O_{\mtt{i}_1}A_{i_2}O_{\mtt{i}_1}\cdots A_{i_1}O_{\mtt{i}_1}\\
        ={}&\frac{1}{a_{\mtt{i}_0}a_{\mtt{i}_1}\cdots a_{\mtt{i}_k}}A_{\mtt{i}}
    \end{aligned}
\end{equation}
The following lemma is immediate.
\begin{lemma}\label{l:norms-match}
    For any $\mtt{i}\in\mathcal{I}^*$ and any subspace $V$ of $\R^2$, we have
    \begin{equation*}
        \frac{\alpha_k(A_{\mtt{i}})}{\|A_{\mtt{i}}|V\|}=\frac{\alpha_k(\overline{A}_{[\mtt{i}]})}{\|\overline{A}_{[\mtt{i}]}|V\|},
    \end{equation*}
    for $k=1,2$.
\end{lemma}
\begin{proof}
    It follows from \cref{e:products} that
    \begin{align*}
        &\alpha_1(\overline{A}_{[\mtt{i}]})=\|\overline{A}_{[\mtt{i}]}\|=\frac{1}{|a_{\mtt{i}_0}a_{\mtt{i}_1}\cdots a_{\mtt{i}_k}|}\|A_{\mtt{i}}\|=\frac{1}{|a_{\mtt{i}_0}a_{\mtt{i}_1}\cdots a_{\mtt{i}_k}|}\alpha_1(A_{\mtt{i}}),\\
        &\alpha_2(\overline{A}_{[\mtt{i}]})=\big\|\overline{A}_{[\mtt{i}]}^{-1}\big\|^{-1}=\frac{1}{|a_{\mtt{i}_0}a_{\mtt{i}_1}\cdots a_{\mtt{i}_k}|}\|A_{\mtt{i}}^{-1}\|^{-1}=\frac{1}{|a_{\mtt{i}_0}a_{\mtt{i}_1}\cdots a_{\mtt{i}_k}|}\alpha_2(A_{\mtt{i}}).
    \end{align*}
    and that for any $v\in\R^2$
    \begin{equation*}
        \|\overline{A}_{[\mtt{i}]}v\|=\frac{1}{|a_{\mtt{i}_0}a_{\mtt{i}_1}\cdots a_{\mtt{i}_k}|}\|A_{\mtt{i}}v\|,
    \end{equation*}
    so the claim follows by the definition of the operator norm.
\end{proof}
Next we show that that the matrices indexed by $[\mtt{i}]$ for $\mtt{i}\in\mathcal{I}^*\setminus\mathcal{I}_e^*$ generate the same semigroup as matrices indexed by $\Lambda^*$.
\begin{lemma}\label{l:same-semigroup}
    For any $\bm{\lambda}\in\Lambda^*$, there exists $\mtt{i}\in\mathcal{I}^*\setminus\mathcal{I}_e^*$ such that
    \begin{equation*}
        A_{[\mtt{i}]}=A_{\bm{\lambda}}
    \end{equation*}
\end{lemma}
\begin{proof}
    Let $\bm{\lambda}=[\mtt{j}_1i_1\mtt{i}_1][\mtt{j}_2i_2\mtt{i}_2]\cdots[\mtt{j}_ki_k\mtt{i}_k]\in\Lambda^*$ be arbitrary, and write $\mtt{i}=\mtt{j}_1i_1\mtt{i}_1\mtt{j}_2\cdots\mtt{i}_{k-1}\mtt{j}_ki_k\mtt{i}_k$.
    Then recalling \cref{e:products},
    \begin{equation*}
        \begin{aligned}
            \overline{A}_{\bm{\lambda}}&=\overline{A}_{[\mtt{j}_1i_1\mtt{i}_1]}\overline{A}_{[\mtt{j}_2i_2\mtt{i}_2]}\cdots \overline{A}_{[\mtt{j}_ki_k\mtt{i}_k]}\\
            &=O_{\mtt{j}_1}A_{i_1}O_{\mtt{i}_1}O_{\mtt{j}_2}A_{i_2}O_{\mtt{i}_2}\cdots O_{\mtt{j}_k}A_{i_k}O_{\mtt{i}_k}\\
            &=O_{\mtt{j}_1}A_{i_1}O_{\mtt{i}_1\mtt{j}_2}A_{i_2}O_{\mtt{i}_2\mtt{j}_3}\cdots O_{\mtt{i}_{k-1}\mtt{j}_k}A_{i_k}O_{\mtt{i}_k}\\
            &=\overline{A}_{[\mtt{i}]}.
        \end{aligned}
    \end{equation*}
    as claimed.
\end{proof}
Now let
\begin{equation*}
    \OA= (\overline{A}_{\lambda})_{\lambda\in\Lambda}.
\end{equation*}
The next proposition shows that $\OA$ is dominated and that the Furstenberg directions of $\A$ are determined by $\OA$.
\begin{proposition}\label{p:dominated-directions}
    If $\A$ is weakly dominated, then $\OA$ is dominated, $Y_F(\A)=Y_F(\OA)$ and $X_F(\A)=X_F(\OA)$.
\end{proposition}
\begin{proof}
    In the proof of \cite[Proposition 2.3]{zbl:1452.37039} the authors show that the strongly invariant multicone for $\A_h$ is also strongly invariant for $\OA$.

    We next show that $Y_F(\A)=Y_F(\OA)$.
    Let $Y\in Y_F(\A)$ and (by definition) find a sequence $c_nA_{\mtt{i}_n}\to A$, with $c_n\in\R$ and $\mtt{i}_n\in\mathcal{I}^*$, where $A$ is a rank one linear map with $\im(A)=Y$.
    We first observe, for all sufficiently large $n$, that $\mtt{i}_n\in \mathcal{I}^*\setminus \mathcal{I}_e^*$.
    Suppose for contradiction that, after passing to a subsequence, $\mtt{i}_{n}\in \mathcal{I}_e^*$ for all $n$.
    Then each $A_{\mtt{i}_n}$ is a constant multiple of an orthogonal matrix, and thus $|\det A_{\mtt{i}_n}|=\|A_{\mtt{i}_n}\|^2$.
    Moreover, since the operator norm and determinant are continuous functions from $\MM_2$ to $\R$,
    \begin{equation*}
        |\det A|=\lim_{n\to\infty}|\det c_nA_{\mtt{i}_n}|=\lim_{n\to\infty}|c_n|^2|\det A_{\mtt{i}_n}|=\lim_{n\to\infty}\|c_nA_{\mtt{i}_n}\|^2=\|A\|^2>0.
    \end{equation*}
    Therefore $A$ has rank two, which is a contradiction.
    Now by recalling \cref{e:products}, we may write
    \begin{equation*}
        \left(c_n\prod_{j=1}^ka_{\mtt{i}_j}\right)\overline{A}_{[\mtt{i}_n]}=c_n A_{\mtt{i}_n}\to A,
    \end{equation*}
    so by definition, $Y\in X_F(\OA)$.

    For the other inclusion, take $Y\in X_F(\OA)$ and again find a sequence $c_nA_{\bm{\lambda}_n}\to A$, with $c_n\in\R$ and $\bm{\lambda}_n\in\Lambda^*$, where $A$ is a rank one linear map with $\im(A)=Y$.
    Apply \cref{l:same-semigroup} to find words $\mtt{i}_n\in\mathcal{I}^*\setminus \mathcal{I}_e^*$ such that $A_{\bm{\lambda}_n}=A_{[\mtt{i}_n]}$ for all $n\in\N$.
    Then by \cref{e:products},
    \begin{equation*}
        \frac{c_n}{\prod_{j=1}^ka_{\mtt{i}_j}}A_{\mtt{i}_n}=c_n \overline{A}_{[\mtt{i}_n]}\to A,
    \end{equation*}
    so $Y\in Y_F(\A)$.

    Finally, if $\A$ is weakly dominated, then so is $\A^{-1}$ and clearly $\overline{\A^{-1}}=\OA^{-1}$.
    Therefore $X_F(\A)=Y_F(\A^{-1})=Y_F(\OA^{-1})=X_F(\OA)$.
\end{proof}
Combining \cref{l:norms-match,l:same-semigroup,p:dominated-directions} gives the following analogue of \cref{l:dominated-restrictions} for weakly dominated matrices.
\begin{lemma}\label{l:weakly-dominated-restrictions}
    Suppose $\A$ is weakly dominated.
    Then there exists a constant $C_1\geq 1$ such that
    \begin{equation*}
        \|A_{\mtt{i}}|Y\|\leq \alpha_1(A_{\mtt{i}})\leq C_1\|A_{\mtt{i}}|Y\|,
    \end{equation*}
    for all $\mtt{i}\in\mathcal{I}^*$ and $Y\in Y_F(\A)$.
    If $V\in X_F(\A)$, then there exists a sequence $\mtt{i}_n\in\mathcal{I}^*$ such that $A_{\mtt{i}_n}V\in X_F(\A)$, and
    \begin{equation*}
        C_1^{-1}\|A_{\mtt{i}_n}|V\|\leq \alpha_2(A_{\mtt{i}_n})\leq \|A_{\mtt{i}_n}|V\|,
    \end{equation*}
    for all $n\in\N$.
    Moreover, we have
    \begin{equation*}
        \lim_{n\to\infty}\frac{\alpha_2(A_{\mtt{i}_n})}{\alpha_1(A_{\mtt{i}_n})}=0.
    \end{equation*}
\end{lemma}
\begin{proof}
    The lower bound in the first claim and the upper bound in the second claim are trivial.
    Let us start by proving the upper bound in the first claim.
    Let $\mtt{i}\in\mathcal{I}^*$.
    If $\mtt{i}\in\mathcal{I}_e^*$, then the claim is trivial since $A_{\mtt{i}}$ is a constant multiple of an orthogonal matrix.
    Therefore we may assume that $\mtt{i}\in\mathcal{I}^*\setminus\mathcal{I}_e^*$.
    Since $\OA$ is dominated, by the first claim in \cref{l:dominated-restrictions}, there exists a constant $C\geq 1$ such that for all $[\mtt{i}]\in\Lambda^*$,
    \begin{equation*}
        1\leq \frac{\alpha_1(\overline{A}_{[\mtt{i}]})}{\|\overline{A}_{[\mtt{i}]}|Y\|}\leq C,
    \end{equation*}
    for all $Y\in Y_F(\OA)$.
    Therefore the first claim follows by \cref{l:norms-match,p:dominated-directions}.

    On the other hand, by the second claim in \cref{l:dominated-restrictions}, for any $V\in X_F(\A)=X_F(\OA)$ we may choose a word $\bm{\lambda}\in\Sigma(\Lambda)$ such that
    \begin{equation*}
        C^{-1}\leq\frac{\alpha_2(\overline{A}_{\overleftarrow{\bm{\lambda}|_n}})}{\|\overline{A}_{\overleftarrow{\bm{\lambda}|_n}}|V\|}\leq1,
    \end{equation*}
    for all $n\in\N$.
    Apply \cref{l:same-semigroup} to find words $\mtt{i}_n\in\mathcal{I}^*$ such that $\overline{A}_{[\mtt{i}_n]}=\overline{A}_{\overleftarrow{\bm{\lambda}|_n}}$ for all $n\in\N$.
    The second claim then follows from \cref{l:dominated}~(\cref{im:shift-inv}), \cref{l:norms-match,p:dominated-directions}.
    Finally, since $\OA$ is dominated,
    \begin{equation*}
        \lim_{n\to\infty}\frac{\alpha_2(A_{\mtt{i}_n})}{\alpha_1(A_{\mtt{i}_n})}=\lim_{n\to\infty}\frac{\alpha_2(\overline{A}_{[\mtt{i}_n]})}{\alpha_1(\overline{A}_{[\mtt{i}_n]})}=\lim_{n\to\infty}\frac{\alpha_2(\overline{A}_{\overleftarrow{\bm{\lambda}|_n}})}{\alpha_1(\overline{A}_{\overleftarrow{\bm{\lambda}|_n}})}=0,
    \end{equation*}
    which is the last claim.
\end{proof}
We recall that a matrix $A\in\MM_2$ has rank one if and only if there are $v\in\im(A)$ and $w\in\ker(A)$ such that $A=vw^{\top}$.
It is then easy to see that
\begin{equation*}
    A=\begin{cases}
        \langle v,w\rangle\pi_{\im(A)}^{\ker(A)},&\text{if $A^2\ne 0$}\\
        \|v\|\|w\|R\pi_{\ker(A)^{\bot}}&\text{if $A^2=0$},
    \end{cases}
\end{equation*}
where $R$ is a rotation by angle $\pi/2$; see, for instance, \cite[Lemma~2.1]{zbmath:7950030}.
By \cite[Lemma 3.2]{zbl:1452.37039}, for weakly dominated tuples $\A$, if $A\in \mathfrak{R}(\A)$ is arbitrary, then $A^2\neq 0$ so in particular $A$ is of the form $\kappa\pi_{\im(A)}^{\ker(A)}$, for some $\kappa\in\R$.
Let us record one final lemma.
\begin{lemma}\label{l:X_F-are-kernels}
    For any $V\in X_F(\A)$ there exists $\kappa\in \R$ and a sequence $\mtt{j}_k\in\mathcal{I}^*$ such that $|\mtt{j}_k|\to \infty$ and
    \begin{equation*}
        \frac{A_{\mtt{j}_k}}{\|A_{\mtt{j}_k}\|}\to \kappa\pi_Y^V,
    \end{equation*}
    in the topology of uniform convergence, where $Y\in Y_F(\A)$.
\end{lemma}
\begin{proof}
    Let $V\in X_F(\A)$ and apply \cref{l:weakly-dominated-restrictions} to find a sequence $\mtt{i}_n\in\mathcal{I}^*$ such that
    \begin{equation*}
        C_1^{-1}\|A_{\mtt{i}_n}|V\|\leq \alpha_2(A_{\mtt{i}_n})\leq \|A_{\mtt{i}_n}|V\|,
    \end{equation*}
    for all $n\in\N$ and
    \begin{equation*}
        \lim_{n\to\infty}\frac{\alpha_2(A_{\mtt{i}_n})}{\alpha_1(A_{\mtt{i}_n})}=0.
    \end{equation*}
    Noting that
    \begin{equation*}
        \left|\det \frac{A_{\mtt{i}_n}}{\|A_{\mtt{i}_n}\|}\right|=\frac{|\det A_{\mtt{i}_n}|}{\|A_{\mtt{i}_n}\|^2}=\frac{\alpha_2(\overline{A}_{[\mtt{i}_n]})}{\alpha_1(\overline{A}_{[\mtt{i}_n]})}\to 0,
    \end{equation*}
    and that $\|A_{\mtt{j}_k}\|^{-1}A_{\mtt{j}_k}\in\{A\in\MM_2\colon \|A\|=1\}$, which is a compact set, by possibly passing to a subsequence, we find a sequence $\mtt{j}_k\coloneqq \mtt{i}_{n_k}$ such that $\|A_{\mtt{j}_k}\|^{-1}A_{\mtt{j}_k}$ converges to a rank one matrix $A$ with $Y\coloneqq \im(A)\in Y_F(\A)$.
    Moreover, for any $v\in V$, we have
    \begin{equation*}
        \frac{\|A_{\mtt{j}_k}v\|}{\|A_{\mtt{j}_k}\|}\leq C_1\frac{\alpha_2(A_{\mtt{j}_k})}{\alpha_2(A_{\mtt{j}_k})}\|v\|\to 0,
    \end{equation*}
    so the kernel of the rank one limit map is $V$ and the claim follows.
\end{proof}

\section{Slicing self-affine sets}\label{s:proof}
In this section, we study planar affine IFSs $(T_i)_{i\in\mathcal{I}}$, that is for every $i\in\mathcal{I}$, $T_i(x)=A_ix+b_i$, for some $A_i\in \GL_2$ and $b_i\in\R^2$.
The attractor of $(T_i)_{i\in\mathcal{I}}$ is denoted by $K$.
We abuse terminology slightly by saying that $K$ is weakly dominated if the tuple $\A=(A_i)_{i\in\mathcal{I}}$ of the linear parts of the associated IFS is weakly dominated.
Since in this section $\A$ is always the tuple of the linear parts of the IFS, we denote the forward and backward Furstenberg directions of $\A$ simply by $X_F$ and $Y_F$, respectively.

The following result is the main goal of this section, and makes up the majority of the proof of \cref{it:tan-slice-form}.
For the reader less familiar with matrix products, we have also included a proof of this result for diagonal systems in \cref{as:diag} which still captures the essence of the pigeonholing argument.
\begin{theorem}\label{t:weakly-dominated-assouad}
    Let $K$ be a weakly dominated self-affine set.
    Then the function $V \mapsto \dimA \pi_{V^\bot}(K)$ is constant on $X_F$.
    Moreover, if $\eta$ denotes this constant value, then for all $V\in X_F$,
    \begin{equation*}
        \dimA K = \eta +\sup_{E\in\Tan(K)}\sup_{x\in\pi_{V^{\bot}}(K)}\dimA(\pi_{V^{\bot}}^{-1}(x)\cap E).
    \end{equation*}
\end{theorem}
We start the proof with a sequence of geometric lemmas.
\subsection{Geometric lemmas}\label{ss:geom}
For the remainder of the paper, for $V\in\RP^1$, we let $e_V$ denote the unique unit vector in $V$ with positive $x$-coordinate if $V$ is not the $y$-axis, and otherwise we let $e_V=(0,1)$.
Going forward, for any $x\in V\in\RP^1$ and $W\in\RP^1$, we let
\begin{equation*}
    x_W=\langle x, e_V\rangle e_W.
\end{equation*}
Similarly, for any subset $E\subset V\in\RP^1$ and $W\in \RP^1$, we let
\begin{equation*}
    E_W=\{x_W\colon x\in E\},
\end{equation*}
that is $E_W$ is a similar copy of $E$ on $W$.
If $x\in V\in\RP^1$, then for any $2\times 2$ matrix $A$, we have
\begin{equation*}
    Ax=\pm\|A|V\|x_{AV},
\end{equation*}
where the sign is positive if $Ae_V$ has the same direction as $e_{AV}$ and negative otherwise.
In our arguments, namely in the proof of \cref{t:weakly-dominated-assouad}, we only need to consider the matrices $A_{\mtt{i}}$ restricted on a direction $Y$ which we fix in the beginning.
Since our arguments rely on pigeonholing and since in any collection of linear maps, at least one half of the maps have a common sign in the equation above, to simplify the situation slightly, we will assume without loss of generality that,
\begin{equation*}
    A_{\mtt{i}}|V(x)=\|A|V\| x_{AV}
\end{equation*}
for all $\mtt{i}\in \mathcal{I}^*$ and $V\in \RP^1$.

Using this assumption, we observe that if a planar set is close to a subset $E$ of a line in $\R^2$, then its image is close to a scaled and rotated copy of $E$.
\begin{lemma}\label{l:image-of-almost-line}
    Let $E\subset Y\in\RP^1$ and let $B\subset \R^2$.
    If $p_{\mathcal{H}}(E;B)<\varepsilon$, then for any $\mtt{i}\in\mathcal{I}^*$,
    \begin{equation*}
        p_{\mathcal{H}}\left(\|A_{\mtt{i}}|Y\|E_{A_{\mtt{i}}Y};A_{\mtt{i}}(B)\right)<\|A_{\mtt{i}}\|\varepsilon.
    \end{equation*}
\end{lemma}
Next we see that when restricted to nearby lines in $Y_F$, all matrices in the semigroup generated by $\A$ have contraction ratios uniformly close to each other.
\begin{lemma}\label{l:norm-close-if-angle-close}
    If $\A$ is weakly dominated, then there exists $C_2>0$ such that the following holds: For all $Y_1,Y_2\in Y_F$ and $\mtt{i}\in\mathcal{I}^*$, if
    \begin{equation*}
        \sin\measuredangle(Y_1,Y_2)\leq C_2\varepsilon,
    \end{equation*}
    for some $0<\varepsilon<1$, then
    \begin{equation*}
        \Big|\|A_{\mtt{i}}|Y_1\|-\|A_{\mtt{i}}|Y_2\|\Big|\leq \varepsilon\|A_{\mtt{i}}|Y_1\|.
    \end{equation*}
\end{lemma}
\begin{proof}
    Assume without loss of generality that $\|A_{\mtt{i}}|Y_1\|\geq \|A_{\mtt{i}}|Y_2\|$.
    Let $C_2=\frac{1}{2C_1}$, where $C_1$ is the constant of \cref{l:weakly-dominated-restrictions}.
    If
    \begin{equation*}
        \sin\measuredangle(Y_1,Y_2)\leq C_2\varepsilon,
    \end{equation*}
    then it follows that $e_{Y_2}=e_{Y_1}+v$, for some $v\in\R^2$, with $\|v\|\leq \frac{\varepsilon}{C_1}$.
    By \cref{l:weakly-dominated-restrictions}
    \begin{equation*}
        \|A|Y_1\|=\|Ae_{Y_1}\|\leq \|Ae_{Y_2}\|+\|Av\|\leq \|Ae_{Y_2}\|+\|A\|\|v\|\leq \|A|Y_2\|+\varepsilon\|A|Y_1\|,
    \end{equation*}
    which gives the claim.
\end{proof}
The following upper bound for the angle between images of lines in $Y_F$ is immediate by combining \cref{l:weakly-dominated-restrictions} and \cite[III Lemma 4.2]{zbl:0572.60001}.
\begin{lemma}\label{l:directions-converge}
    If $\A$ is weakly dominated, then there is a constant $C_3>0$ such that for all $V,W\in Y_F$ and $\mtt{i}\in\mathcal{I}^*$,
    \begin{equation*}
        \measuredangle(A_{\mtt{i}}V,A_{\mtt{i}}W)\leq C_3\frac{\alpha_2(A_{\mtt{i}})}{\alpha_1(A_{\mtt{i}})}
    \end{equation*}
\end{lemma}
The next lemma is a trivial consequence of the fact that there is a unique way to decompose any $x\in \R^2$ as the sum of vectors in $V,W\in\RP^1$ with $V\ne W$ given by
\begin{equation*}
    x=\pi_{V}^{W}(x)+\pi_{W}^{V}(x).
\end{equation*}
\begin{lemma}\label{l:diam-of-projection}
    For all $x\in\R^2$, $V,W\in \RP^1$, with $V\ne W$ and $\mtt{i}\in\mathcal{I}^*$, we have
    \begin{equation*}
        \pi_{A_{\mtt{i}}V}^{A_{\mtt{i}}W}(A_{\mtt{i}}x)= \|A_{\mtt{i}}|V\|\pi_{V}^{W}(x)_{A_{\mtt{i}}V}.
    \end{equation*}
\end{lemma}
Finally, an elementary geometric argument shows that projections along the directions in $Y_F$ of well separated sets on lines in the directions of $X_F$ are well separated.
\begin{lemma}\label{l:separated-projections}
    Let $A\subset V\in X_F$ be a $(1+\frac{2}{\sin\delta})r$-separated set, where $\delta=\measuredangle(X_F,Y_F)>0$, and for each $a\in A$, let $\xi(a)\in B(a,2r)\subset \R^2$.
    Then for any $Y\in Y_F$, the set
    \begin{equation*}
        \pi_V^Y(\{\xi(a)\colon a\in A\}),
    \end{equation*}
    is $r$-separated.
\end{lemma}
\begin{proof}
    Since $A\subset V$, $\pi_V^Y(A)=A$.
    Moreover, for any $Y\in Y_F$ and $a\in A$, the set $\pi_V^Y(B(a,2r))$ is an interval of width $2(\sin\measuredangle(V,Y))^{-1} r\leq2(\sin\delta)^{-1} r$ centred at $a$.
    Therefore,
    \begin{equation*}
        |\xi(a)-\xi(b)|\geq \left(1+\frac{2}{\sin\delta}\right)r-\frac{2}{\sin\delta}r=r,
    \end{equation*}
    for any $a,b\in A$.
\end{proof}

\subsection{Product structure of weak tangents}
We are now ready to prove \cref{t:weakly-dominated-assouad}.
The main part of the proof is in \cref{t:product-tangent}, which gives a slightly stronger result than just the lower bound in \cref{t:weakly-dominated-assouad}, namely that the Assouad dimension is attained by a coarse microset which is a product set.
Before proving this theorem, let us observe that the map $V\mapsto \dimA\pi_{V^{\bot}}(K)$ is constant on $X_F$.
\begin{proposition}\label{p:projection-assouad-constant}
    Let $K$ be a weakly dominated self-affine set.
    Then
    \begin{equation*}
        \dimA \pi_{W^{\bot}}(K)=\dimA \pi_{V^{\bot}}(K)
    \end{equation*}
    for all $V,W\in X_F$.
\end{proposition}
\begin{proof}
    Let $V,W\in X_F$.
    It suffices to show that
    \begin{equation*}
        \dimA \pi_{W^{\bot}}(K)\geq \dimA \pi_{V^{\bot}}(K).
    \end{equation*}
    By \cref{l:X_F-are-kernels}, there exists a sequence $\mtt{j}_k\in\mathcal{I}^*$ such that $|\mtt{j}_k|\to \infty$, and
    \begin{equation*}
        \frac{A_{\mtt{j}_k}}{\|A_{\mtt{j}_k}\|}\to \kappa\pi_Y^V,
    \end{equation*}
    in the topology of uniform convergence, where $Y\in Y_F$.
    Note that $\kappa\pi_Y^V(K)$ is bi-Lipschitz equivalent to $\pi_{V^{\bot}}(K)$, so in particular $\dimA \kappa\pi_Y^V(K)=\dimA \pi_{V^{\bot}}(K)$.
    Thus for any $s<\dimA \pi_{V^{\bot}}(K)$ and $C>0$, we may choose a point $x\in \kappa\pi_Y^V(K)$ and a collection $\mathcal{A}\subset \kappa\pi_Y^V(K)\cap B(x,R)\subset Y$ of $3r$-separated points with
    \begin{equation*}
        \#\mathcal{A}\geq C\left(\frac{R}{r}\right)^s.
    \end{equation*}
    Let $\delta=\measuredangle(Y,W)>0$ and take $k$ sufficiently large so that
    \begin{equation*}
        \left\|\frac{A_{\mtt{j}_k} x}{\|A_{\mtt{j}_k}\|}-\kappa\pi_Y^V(x)\right\|<r\sin{\delta}
    \end{equation*}
    for all $x\in\R^2$.
    Therefore, by translating $K$ if necessary we find for each $a\in \mathcal{A}$ a point $x(a)\in K$ such that
    \begin{equation*}
        \left\|x(a)- a\|A_{\mtt{j}_k}\|\right\|<r\|A_{\mtt{j}_k}\| \sin{\delta}.
    \end{equation*}
    Of course, $\pi_{W^{\bot}}(\mathcal{A})$ is $3r\sin\delta$-separated and therefore the points
    \begin{equation*}
        \mathcal{B}\coloneqq \{\pi_{W^{\bot}}(x(a))\colon a\in\mathcal{A}\},
    \end{equation*}
    form a $r\|A_{\mtt{j}_k}\| \sin\delta $-separated subset of $\pi_{W^{\bot}}(K)\cap B(\pi_{W^{\bot}}(x),R\|A_{\mtt{j}_k}\|\sin\delta )$ with
    \begin{equation*}
        \#\mathcal{B}=\#\mathcal{A}\geq C\left(\frac{R\|A_{\mtt{j}_k}\|}{r\|A_{\mtt{j}_k}\|}\right)^s.
    \end{equation*}
    Since this holds for all $C>0$, by the definition of the Assouad dimension, we have
    \begin{equation*}
        \dimA\pi_{W^{\bot}}(K)\geq s
    \end{equation*}
    and since $s< \dimA \pi_{V^{\bot}}(K)$ was arbitrary, we are done.
\end{proof}
\begin{remark}
    It is possible to calculate the constant value $\eta=\dimA\pi_{V^{\bot}}(K)$ in the proposition above in some situations.
    If $K$ is a self-affine carpet, then the set $X_F$ is a singleton consisting of the direction with the strongest contraction ratio.
    In this case, the projection of the IFS along this direction is a self-similar IFS on the line and therefore, by \cite{zbl:1317.28014}, $\eta=\dimH \pi(K)$ if the projected IFS satisfies the \emph{weak separation condition} and $\eta=1$ otherwise.

    On the other hand, for general self-affine sets, it is easy to see by using  \cite[Lemma 6.4]{arxiv:2107.00983} that in the absence of a projective separation condition, which is similar in spirit to the open set condition, we have that $\eta=1$.
    Moreover this projective separation condition is generically, in a topological sense, not satisfied by self-affine sets satisfying the strong separation condition.
    Formally, it was shown in \cite[Theorem 3.6]{arxiv:2107.00983} that for a given collection of linear parts $(A_i)_{i\in\mathcal{I}}$, there is a \emph{residual set}---a countable intersection of sets with dense interiors---of translation vectors $(t_i)_{i\in\mathcal{I}}$, such that the IFS $(x\mapsto A_ix+t_i)_{i\in\mathcal{I}}$ does not satisfy the projective open set condition.
    Characterizing the projections in a similar fashion to the carpet setting would be of interest but we do not pursue this further in this work.
\end{remark}

\begin{figure}[t]
    \begin{subcaptionblock}{0.47\textwidth}
        \centering
        \begin{tikzpicture}[scale=0.7]
    \clip (-1,0) rectangle + (10,12);
    \pgfmathsetseed{5}
    \coordinate (start) at (2,0);
    \coordinate (left) at (1.5,0);
    \coordinate (right) at (2.5,0);
    \def\rangle{85} %

    \draw[dashed] (left) -- ++(\rangle:10.4) node[above right] {$V\in X_F$};
    \draw[dashed] (right) -- ++(\rangle:10.4);

    \foreach \i in {2,3,4,5,6,10,11,12,17,18} {
        \pgfmathsetmacro{\ellipseAngle}{int(rand*25)}
        \pgfmathsetmacro{\offsetDist}{rand*0.6}
        \path (start) ++(\rangle:{\i * 0.5}) ++(\ellipseAngle:0.5*\offsetDist) coordinate (P\i);
        \path (P\i) ++(\ellipseAngle:\offsetDist) coordinate (E\i);
        \draw[rotate around={\ellipseAngle:(E\i)}, color=gray] (E\i) ellipse (0.6 and 0.1);
        \draw[color=gray] (E\i) -- ++(\ellipseAngle:0.3);
    }
    \foreach \i in {1,7,8,9,13,14,15,19} {
        \pgfmathsetmacro{\ellipseAngle}{int(rand*25)}
        \pgfmathsetmacro{\offsetDist}{rand*0.6}
        \path (start) ++(\rangle:{\i * 0.5}) ++(\ellipseAngle:0.5*\offsetDist) coordinate (P\i);
        \path (P\i) ++(\ellipseAngle:\offsetDist) coordinate (E\i);
        \draw[rotate around={\ellipseAngle:(E\i)}, color=gray] (E\i) ellipse (0.6 and 0.2);
        \draw[color=gray] (E\i) -- ++(\ellipseAngle:0.3);
    }
    \foreach \i in {2,6,8,9} {
        \path (start) ++(\rangle:{\i}) coordinate (P\i);
        \pgfmathsetmacro{\ellipseAngle}{20}
        \pgfmathsetmacro{\offsetDist}{0.2}
        \path (P\i) ++(\ellipseAngle:\offsetDist) coordinate (E\i);
        \draw[rotate around={\ellipseAngle:(E\i)}] (E\i) ellipse (0.6 and 0.1);
        \draw[thick] (E\i) -- ++(\ellipseAngle:0.3);
    }
    \foreach \i in {1,3} {
        \path (start) ++(\rangle:{\i}) coordinate (P\i);
        \pgfmathsetmacro{\ellipseAngle}{20}
        \pgfmathsetmacro{\offsetDist}{0.2}
        \path (P\i) ++(\ellipseAngle:\offsetDist) coordinate (E\i);
        \draw[rotate around={\ellipseAngle:(E\i)}] (E\i) ellipse (0.6 and 0.2);
        \draw[thick] (E\i) -- ++(\ellipseAngle:0.3);
    }
\end{tikzpicture}
        \caption{The initial pigeonholed cylinders.}
        \label{f:main-figure-1}
    \end{subcaptionblock}
    \begin{subcaptionblock}{0.47\textwidth}
        \centering
        \begin{tikzpicture}[scale=0.7]
    \clip (-2,0) rectangle + (10,12);
    \pgfmathsetseed{5}
    \coordinate (start) at (0,0);
    \def\langle{100} %

    \draw[dashed] (start) ++ (4:1.6) -- ++(\langle:10.4) node[above right] {$A_{\mtt{i}}V$};

    \foreach \i in {1,2,3,5,6,7,9,10,11,13,14,15,17,18,19} {
        \pgfmathsetmacro{\ellipseAngle}{int(rand*25)}
        \pgfmathsetmacro{\offsetDist}{rand*0.6}

        \path (start) ++(\langle:{\i * 0.5}) ++(0.2*\ellipseAngle:3.5*\offsetDist) coordinate (P\i);
        \path (P\i) ++(0.2*\ellipseAngle:7*\offsetDist) coordinate (E\i);
        \draw[rotate around={0.2*\ellipseAngle:(E\i)}, color=gray] (E\i) ellipse (6 and 0.1);
        \draw[color=gray] (E\i) -- ++(0.2*\ellipseAngle:5);
    }
    \foreach \i in {1,...,5} {
        \pgfmathsetmacro{\ellipseAngle}{20}
        \pgfmathsetmacro{\offsetDist}{0.2}

        \path (start) ++(\langle:{\i * 2}) coordinate (P\i);
        \path (P\i) ++(0.2*\ellipseAngle:8*\offsetDist) coordinate (E\i);
        \draw[rotate around={0.2*\ellipseAngle:(E\i)}] (E\i) ellipse (6 and 0.1);
        \draw[thick] (E\i) -- ++(0.2*\ellipseAngle:5);
    }
\end{tikzpicture}
        \caption{The image of the configuration.}
        \label{f:main-figure-2}
    \end{subcaptionblock}
    \caption{Depiction of the proof of \cref{t:product-tangent}.
        In \cref{f:main-figure-1}, we see the initial pigeonholed collection of cylinders (and large pieces of the projections) corresponding to an approximation of a slice of a large weak tangent.
        In \cref{f:main-figure-2}, we see the image of this configuration under an appropriately chosen affine map after locating an improved scale with \cref{l:h-micro}.
    }
    \label{f:main-figure}
\end{figure}
We now move on to the proof of our main result.
Let us begin with an informal overview of the strategy, see \cref{f:main-figure} for an illustration.
We begin by discretizing a given slice of a weak tangent in a backward Furstenberg direction, and attach to each point in the discretized slice a copy of (an approximation of) a weak tangent of the projection using the self-affine structure.
These copies are depicted by the short lines in the ellipses in \cref{f:main-figure-1}.
Unfortunately, the individual approximations of weak tangents of the projection need not line up at all, and since the Assouad dimension can be substantially larger than the Hausdorff dimension, this will cause loss in dimension.
By pigeonholing (with constants depending on the resolution of approximation of the large weak tangent of the projection), we can find a sub-family which has substantially improved alignment properties (the dark ellipses in \cref{f:main-figure-1}).
However, the corresponding component of the slice corresponding to this pigeonholed sub-family could be substantially smaller at the original scale.
In order to amplify this configuration, we apply the discretized version Furstenberg's construction of large microsets (that is, \cref{l:h-micro}) in the slice to locate a new scale where there is no loss in dimension.
Here, it is important that the initial configuration is a very tall and narrow tube, with eccentricity sufficiently large depending on the constants in \cref{l:h-micro} and the resolution of the weak tangent in the projection so that the pigeonholing and amplification step do not flatten the tube beyond being a square.
Of course, since \cref{l:h-micro} does not guarantee a precise scale at which the improved configuration appears, this amplification process could still result in a very thin tube.
Finally, we use the fact that the slice is in a backwards Furstenberg direction and apply an appropriate high iteration $T_{\mtt{i}}$ of the affine maps to squash this tube to a square, see \cref{f:main-figure-2}.
Passing to the limit yields a weak tangent with the desired product structure.
\begin{theorem}\label{t:product-tangent}
    Let $K$ be a weakly dominated self-affine set, $V,W\in X_F$, $E\in \Tan(K)$ and $x\in\pi_{V^{\bot}}(E)$.
    Then there are compact sets $A,B\in\R$ with $\dimlB A=\dimA \pi_{W^{\bot}}(K)$ and $\dimlB B=\dimA(\pi_{V^{\bot}}^{-1}(x)\cap E)$ such that $A\times B$ is a coarse microset of $K$.
\end{theorem}
\begin{proof}
    Let $V,W\in X_F$, $E\in\Tan(K)$ and $x\in\pi_{V^{\bot}}(E)$.
    Write
    \begin{equation*}
        \eta=\dimA \pi_{W^{\bot}}(K)\qquad\text{and}\qquad\beta=\dimA(\pi_{V^{\bot}}^{-1}(x)\cap E).
    \end{equation*}
    First, in the same way as shown in \cref{p:projection-assouad-constant}, we may use \cref{l:X_F-are-kernels}, to find a sequence $\mtt{j}_k\in\mathcal{I}^*$ such that $|\mtt{j}_k|\to \infty$ and
    \begin{equation*}
        \lim_{k\to\infty}\frac{A_{\mtt{j}_k}}{\|A_{\mtt{j}_k}\|} = \kappa\pi_Y^W
    \end{equation*}
    for some $Y\in Y_F$.
    Moreover, recall that
    \begin{equation*}
        \dimA \kappa\pi_Y^W(K)=\dimA \pi_{W^{\bot}}(K) = \eta.
    \end{equation*}
    To simplify notation, let us assume that $\kappa=1$.
    We begin start by constructing approximate copies of the thickest parts of $\pi_{Y}^W(K)$ inside the self-affine set $K$.
    Let $m\in\N$.
    Going forward, most of the choices we make depend on $m$, but to simplify notation this dependence is often left implicit.

    First, by \cref{c:weak-tan}, choose $k\in\N$, scales $0<r=2^{-(k+m)}<R=2^{-k}$ and a point $x\in \pi_{Y}^W(K)$ such that $\mathcal{P}_m\coloneqq \pi_{Y}^W(K)\cap B(x,R)$ satisfies
    \begin{equation}\label{e:proj-packing}
        N_{k+n}(\mathcal{P}_m)\geq \left(\frac{R}{2^{-(k+n)}}\right)^{\eta-\frac{1}{m}}=2^{n(\eta-\frac{1}{m})},
    \end{equation}
    for all $0\leq n\leq m$.
    Now choose $k\in\N$ large enough such that the word $\mtt{j}\coloneqq \mtt{j}_k$ satisfies
    \begin{equation*}
        \left\|\frac{A_{\mtt{j}}x}{\|A_{\mtt{j}}\|}- \pi_{Y}^W x\right\|<\frac{1}{m}R,
    \end{equation*}
    for all $x\in\R^2$.
    Note that $T_{\mtt{j}}(K)-T_{\mtt{j}}(x)=A_{\mtt{j}}(K-x)$, so the previous equation implies that
    \begin{equation}\label{e:almost-line}
        p_{\mathcal{H}}\left(\frac{\mathcal{P}_m-x}{R};\frac{T_{\mtt{j}}(K)-T_{\mtt{j}}(x)}{\|A_{\mtt{j}}\|R}\right)\leq \frac{1}{m}.
    \end{equation}

    Next, set
    \begin{equation*}
        M\geq \frac{8\pi m^3(1-C_1^{-1}C^{-1}\alpha_{\min})}{C_2\|A_{\mtt{j}}\|R^2\sin\delta}
    \end{equation*}
    where $\delta=\measuredangle(X_F,Y_F)$ and $\alpha_{\min}=\min_{i\in\mathcal{I}}\alpha_1(A_i)$.
    Using the definition of the Assouad dimension, find scales $0<r_0<R_0<1$, with $R_0\geq2^{m(m+k)}r_0$, a point $z_0\in \pi_{V^{\bot}}^{-1}(x)\cap E$ and a finite set $\mathcal{A}_0\subset \pi_{V^{\bot}}^{-1}(x)\cap E\cap B(z_0,R_0)$ of $(1+\frac{2}{\sin\delta})r_0$-separated points, satisfying
    \begin{equation*}
        \#\mathcal{A}_0\geq 4M\left(\frac{R_0}{r_0}\right)^{\beta-\frac{1}{m}}.
    \end{equation*}
    Since $E\in\Tan(K)$, there exists $y\in K$ and $\lambda\geq1$ so that
    \begin{equation*}
        d_{\mathcal{H}}(\lambda(K-y)\cap B(0,1),E)\leq \frac{r_0}{2}.
    \end{equation*}
    Choose $\ell_{m,1}\in\N$ maximal and $\ell_{m,2}\in\N$ minimal, so that
    \begin{equation*}
        2^{-\ell_{m,2}}\leq \lambda^{-1}r_0\leq \lambda^{-1}R_0\leq 2^{-\ell_{m,1}},
    \end{equation*}
    and let $r_1=2^{-\ell_{m,2}}$ and $R_1=2^{-\ell_{m,1}}$.
    Now, since the points in $\mathcal{A}_0$ are $(1+\frac{2}{\sin\delta})r_0$-separated, by the choice of $\ell_{m,2}$, the points in $\lambda^{-1}\mathcal{A}_0+y$ are $(1+\frac{2}{\sin\delta})r_1$-separated.
    Moreover, for each $a\in\lambda^{-1}\mathcal{A}_0+y$, there exists $\xi(a)\in K\cap B(a,r_1)$.
    Let us denote by $\mathcal{A}\coloneqq\{\xi(a)\colon a\in\lambda^{-1}\mathcal{A}_0+y\}$.

    For each $y\in\mathcal{A}$, let $\mtt{i}_y\in\mathcal{I}^*$ satisfy $\alpha_1(A_{\mtt{i}_y})\leq r_1 < \alpha_1(A_{\mtt{i}_y^-})$ and $y\in T_{\mtt{i}_y}(K)$.
    To simplify notation, write $x_y=T_{\mtt{i}_y\mtt{j}}(x)$.
    Since each $T_{\mtt{i}_y}(K)$ has diameter at most $r_1$ and the set $\mathcal{A}$ is a subset of a $2r_1$-neighbourhood of a line in direction $V$, we have $\diam(\pi_{Y}^{V}(\bigcup_{y\in\mathcal{A}}T_{\mtt{i}_y}(K)))\leq \frac{4r_1}{\sin\delta}$.
    Additionally, $\diam(\RP^1)=\frac{\pi}{2}$ and $C_1^{-1}C^{-1}\alpha_{\min}r_1\leq \|A_{\mtt{i}_y}|Y\|\leq r_1$, where $C$ and $C_1$ are the constants of \cref{l:almost-multiplicative} and \cref{l:weakly-dominated-restrictions}, respectively, so by the pigeonhole principle there is a line $Y_m\in Y_F$, a real number $w_m\in[C_1^{-1}C^{-1}\alpha_{\min},1]$ and a subset $\mathcal{B}_0\subset \mathcal{A}$, satisfying
    \begin{equation*}
        \#\mathcal{B}_0\geq \frac{1}{M}\cdot\#\mathcal{A}\geq 4\left(\frac{\lambda^{-1}R_0}{\lambda^{-1}r_0}\right)^{\beta-\frac{1}{m}}\geq\left(\frac{R_1}{r_1}\right)^{\beta-\frac{1}{m}},
    \end{equation*}
    such that the following inequalities hold for all $y,z\in\mathcal{B}_0$:
    \begin{align}
        \sin\measuredangle(A_{\mtt{i}_y}Y,Y_m)&\leq C_2\frac{1}{2m},\label{e:dir-close}\\
        \left\|\pi_{Y_m}^{V}(x_y-x_z)\right\|&\leq \frac{\|A_{\mtt{j}}\|Rr_1}{m},\label{e:proj-close}\\
        \Big|\|A_{\mtt{i}_y}|Y\|-w_mr_1\Big|&\leq \frac{Rr_1}{2m}\label{e:width-close}.
    \end{align}
    Here, $C_2>0$ is the constant of \cref{l:norm-close-if-angle-close}.
    Since $\diam(T_{\mtt{i}_{\xi(a)}}(K))\leq r_1$ for any $a\in\lambda^{-1}\mathcal{A}_0+y$, we have
    \begin{equation*}
        x_{\xi(a)}\in B(a,2r_1),
    \end{equation*}
    so by applying \cref{l:separated-projections}, we see that the set
    \begin{equation*}
        \mathcal{B}\coloneqq \{\pi_{V}^{Y_m}(x_y)\colon y\in\mathcal{B}_0\},
    \end{equation*}
    is an $r_1$-separated subset of $V$.

    Now, since the choice of $r_1$ and $r_2$ implies that $\ell_{m,2}-\ell_{m,1}\geq m(m+k)$, we may apply \cref{l:h-micro} with $s=\beta-\frac{2}{m},t=\beta-\frac{1}{m},\ell=m$ and $k=m+k$ to obtain a subcollection $\mathcal{B}_m$ of $\mathcal{B}$ with $\mathcal{B}_m\subset B(\pi_{V}^{Y_m}(x_z),R')$ for some $z\in\mathcal{B}_0$ and $R'=2^{-\ell_{m}}$ with $\ell_{m,1}\leq \ell_{m}\leq \ell_{m,2}-m-k$ such that
    \begin{equation}\label{e:slice-packing}
        N_{\ell_m+n}(\mathcal{B}_m)\geq 2^{n(\beta-\frac{2}{m})},
    \end{equation}
    for all $0\leq n\leq m$.

    Next we will construct the coarse microset.
    By \cref{l:weakly-dominated-restrictions}, we may choose a sequence $\mtt{i}_j\in\mathcal{I}^*$ such that
    \begin{equation*}
        C_1^{-1}\|A_{\mtt{i}_j}|V\|\leq \alpha_2(A_{\mtt{i}_j})\leq \|A_{\mtt{i}_j}|V\|,
    \end{equation*}
    for all $j$ and
    \begin{equation*}
        \frac{\|A_{\mtt{i}_j}|V\|}{\|A_{\mtt{i}_j}|Y_m\|}\leq C_1^{2}\frac{\alpha_2(A_{\mtt{i}_j})}{\alpha_1(A_{\mtt{i}_j})}\to0.
    \end{equation*}
    Therefore, we may choose $j$ to be the smallest number satisfying
    \begin{equation*}
        \frac{\|A_{\mtt{i}_j}|V\|}{\|A_{\mtt{i}_j}|Y_m\|}\leq \|A_{\mtt{j}}\|R\frac{r_1}{R'},
    \end{equation*}
    and note that then by \cref{l:weakly-dominated-restrictions,l:almost-multiplicative}
    \begin{equation*}
        \frac{\|A_{\mtt{i}_j}|V\|}{\|A_{\mtt{i}_j}|Y_m\|}\geq (C_1^{-1}C^{-1}\alpha_{\min})^2\frac{\|A_{\mtt{i}|_{j-1}}|V\|}{\|A_{\mtt{i}|_{j-1}}|Y_m\|}\geq (C_1^{-1}C^{-1}\alpha_{\min})^2\|A_{\mtt{j}}\|R\frac{r_1}{R'}.
    \end{equation*}
    Therefore, by possibly passing to a subsequence, we may assume that
    \begin{equation*}
        h_m\coloneqq \frac{\|A_{\mtt{i}_j}|V\|R'}{\|A_{\mtt{j}}\|R\|A_{\mtt{i}_j}|Y_m\|r_1}\to h,
    \end{equation*}
    with $h\in [(C_1^{-1}C^{-1}\alpha_{\min})^2,1]$.
    Let us denote $Y_m'=A_{\mtt{i}_j}Y_m$ and $Y_{y,m}=A_{\mtt{i}_j\mtt{i}_y}Y$.
    \cref{l:directions-converge} implies that for all $y\in\mathcal{A}$,
    \begin{equation}\label{e:angle-close}
        \measuredangle(Y_m',Y_{y,m})\leq C\frac{\alpha_2(A_{\mtt{i}_j})}{\alpha_1(A_{\mtt{i}_j})}\leq C\frac{\|A_{\mtt{i}_j}|V\|}{\|A_{\mtt{i}_j}|Y_m\|}\leq C\|A_{\mtt{j}}\|R\frac{r_1}{R'}\leq C\|A_{\mtt{j}}\|2^{-m}R.
    \end{equation}
    In particular, we certainly have $\sin\measuredangle(Y_m',Y_{y,m})\leq \frac{1}{m}$, for all large enough $m$.
    Moreover, by \cref{l:weakly-dominated-restrictions}, $V_m\coloneqq A_{\mtt{i}_j}V\in X_F$ and therefore, passing to a subsequence if necessary, we may assume that there are $V'\in X_F$ and $Y'\in Y_F$ such that $V_m\to V'$ and $Y_m'\to Y'$.
    Now \cref{e:proj-close} together with \cref{l:diam-of-projection} shows that
    \begin{equation*}
        \left\|\frac{\pi_{Y_m'}^{V_m}A_{\mtt{i}_j}(x_y-x_z)}{\|A_{\mtt{j}}\|R\|A_{\mtt{i}_j}|Y_m\|r_1}\right\|=\frac{\|A_{\mtt{i}_j}|Y_m\|\|\pi_{Y_m}^{V}(x_y-x_z)\|}{\|A_{\mtt{j}}\|R\|A_{\mtt{i}_j}|Y_m\|r_1}\leq \frac{1}{m}.
    \end{equation*}
    Since $x=\pi_{Y_{m}'}^{V_m}(x)+\pi_{V_m}^{Y_{m}'}(x)$ for any $x\in\R^2$, the previous equation gives
    \begin{equation}\label{e:point-close}
        \begin{aligned}
            \MoveEqLeft\left\|\frac{T_{\mtt{i}_j}(x_y)-T_{\mtt{i}_j}(x_z)}{\|A_{\mtt{j}}\|R\|A_{\mtt{i}_j}|Y_m\|r_1}-h_m\frac{\pi_{V}^{Y_m}(x_y-x_z)_{V_m}}{R'}\right\|\\*
            &=\left\|\frac{A_{\mtt{i}_j}(x_y-x_z)-\|A_{\mtt{i}_j}|V\|\pi_{V}^{Y_m}(x_y-x_z)_{V_m}}{\|A_{\mtt{j}}\|R\|A_{\mtt{i}_j}|Y_m\|r_1}\right\|\\
            &\leq \left\|\frac{\pi_{Y_m'}^{V_m}A_{\mtt{i}_j}(x_y-x_z)}{\|A_{\mtt{j}}\|R\|A_{\mtt{i}_j}|Y_m\|r_1}\right\|+\left\|\frac{\pi_{V_m}^{Y_m'}A_{\mtt{i}_j}(x_y-x_z)-\pi_{V_m}^{Y_m'}A_{\mtt{i}_j}(x_y-x_z)}{\|A_{\mtt{j}}\|R\|A_{\mtt{i}_j}|Y_m\|r_1}\right\|\\
            &\leq \frac{1}{m}.
        \end{aligned}
    \end{equation}
    Since for any $A,B\in \GL_2$, and $Y\in\RP^1$ we have $\|AB|Y\|=\|A|BY\|\|B|Y\|$, by \cref{e:width-close} and \cref{l:norm-close-if-angle-close} together with \cref{e:dir-close},
    \begin{equation}\label{e:widths-close}
        \begin{aligned}
            \MoveEqLeft\Big|\|A_{\mtt{i}_j\mtt{i}_y}|Y\|-w_m\|A_{\mtt{i}_j}|Y_m\|r_1\Big|= \Big|\|A_{\mtt{i}_j}|Y_{y,m}\|\|A_{\mtt{i}|_y}|Y\|-w_m\|A_{\mtt{i}_j}|Y_m\|r_1\Big|\\*
            &\leq\Big|\|A_{\mtt{i}_j}|Y_{y,m}\|-\|A_{\mtt{i}_j}|Y_m\|\Big|\|A_{\mtt{i}|_{y}}|Y\|+\Big|\|A_{\mtt{i}|_y}|Y\|-w_mr_1\Big|\|A_{\mtt{i}_j}|Y_m\|\\
            &\leq \frac{1}{m}\|A_{\mtt{i}_j}|Y_m\|r_1
        \end{aligned}
    \end{equation}
    and again by passing to a subsequence if necessary, we may assume that $w_m\to w$ for some $w\in[C_1^{-1}C^{-1}\alpha_{\min},1]$.

    Let now
    \begin{equation*}
        A_m = \frac{(\mathcal{P}_m-x)_{\R}}{R}\quad\text{ and }\quad B_m = \frac{(\mathcal{B}_m-\pi_{V}^{Y_m}(x_z))_{\R}}{R'}
    \end{equation*}
    and recall that $\mathcal{P}_m\subset B(x,R)$ and $\mathcal{B}_m\subset B(\pi_{V}^{Y_m}(x_z),R')$, so $A_m\times B_m$ is a compact subset of $B(0,1)$.
    Therefore, passing to a subsequence, get a compact set $A\times B$ such that $A_m\times B_m\to A\times B$ in the Hausdorff distance.
    By \cref{e:proj-packing}, \cref{e:slice-packing} and \cref{l:lower-box}, we see that $\dimlB A\geq \eta$ and $\dimlB B\geq \beta$.

    Finally, let us show that $A\times B$ is a coarse microset of $K$.
    Let $f_m\colon \R^2\to\R^2$ be the unique linear map taking the vector $(1,0)$ to $w_me_{Y_m'}$ and $(0,1)$ to $h_me_{V_m}$, and let $f\colon \R^2\to\R^2$ be the unique linear map taking the vector $(1,0)$ to $we_{Y'}$ and $(0,1)$ to $he_{V'}$.
    Clearly the function $f$ is bi-Lipschitz, $f_m\to f$ in the topology of uniform convergence and $f_m(A_m\times B_m)\to f(A\times B)$ in the Hausdorff distance.
    Moreover,
    \begin{equation*}
        f_m(A_m\times B_m)=w_m\frac{(\mathcal{P}_m-x)_{Y_m'}}{R}+h_m\frac{(\mathcal{B}_m-\pi_{V}^{Y_m}(x_z))_{V_m}}{R'}.
    \end{equation*}
    Therefore, to finish the proof it suffices to show that
    \begin{equation}\label{e:weak-pseudo-tangent}
        \lim_{m\to\infty} p_{\mathcal{H}}\left(f_m(A_m\times B_m);\frac{T_{\mtt{i}_{j}}(K)-T_{\mtt{i}_{j}}(x_z)}{\|A_{\mtt{j}_n}\|R\|A_{\mtt{i}_{j}}|Y_m\|r_1}\right)=0.
    \end{equation}
    Note that for each $\pi_{V}^{Y_m}(x_y)\in\mathcal{B}_m$,
    \begin{equation*}
        f_m\left(\frac{(\mathcal{P}_m-x)_{Y_m'}}{R}+\frac{\pi_{V}^{Y_m}(x_y-x_z)_{V_m}}{R'}\right)=w_m\frac{(\mathcal{P}_m-x)_{Y_m'}}{R}+h_m\frac{\pi_{V}^{Y_m}(x_y-x_z)_{V_m}}{R'}.
    \end{equation*}
    By combining \cref{e:almost-line} with \cref{l:image-of-almost-line}, and recalling that $x_y=T_{\mtt{i}_y\mtt{j}}(x)$ and $w_m\leq 1$,
    \begin{equation*}
        p_{\mathcal{H}}\left(w_m\frac{(\mathcal{P}_m-x)_{Y_{y,m}}}{R};w_m\frac{T_{\mtt{i}_{j}\mtt{i}_y\mtt{j}}(K)-T_{\mtt{i}_{j}}(x_y)}{\|A_{\mtt{j}}\|R\|A_{\mtt{i}_{j}\mtt{i}_y}|Y\|}\right)\leq \frac{C_1}{m},
    \end{equation*}
    and by \cref{e:widths-close},
    \begin{equation*}
        d_{\mathcal{H}}\left(w_m\frac{T_{\mtt{i}_{j}\mtt{i}_y\mtt{j}}(K)-T_{\mtt{i}_{j}}(x_y)}{\|A_{\mtt{j}}\|R\|A_{\mtt{i}_{j}\mtt{i}_y}|Y\|},\frac{T_{\mtt{i}_{j}\mtt{i}_y\mtt{j}}(K)-T_{\mtt{i}_{j}}(x_y)}{\|A_{\mtt{j}}\|R\|A_{\mtt{i}_{j}}|Y_m\|r_1}\right)\leq \frac{1}{m}.
    \end{equation*}
    Since $\sin\measuredangle(Y_m',Y_{y,m})\leq \frac{1}{m}$ and $\mathcal{P}_m\in B(x,R)$,
    \begin{equation*}
        d_{\mathcal{H}}\left(w_m\frac{(\mathcal{P}_m-x)_{Y_m'}}{R},w_m\frac{(\mathcal{P}_m-x)_{Y_{y,m}}}{R}\right)\leq \frac{1}{m}
    \end{equation*}
    which together with the previous two inequalities imply that
    \begin{equation*}
        p_{\mathcal{H}}\left(w_m\frac{(\mathcal{P}_m-x)_{Y_m'}}{R};\frac{T_{\mtt{i}_{j}\mtt{i}_y\mtt{j}}(K)-T_{\mtt{i}_{j}}(x_y)}{\|A_{\mtt{j}}\|R\|A_{\mtt{i}_{j}}|Y_m\|r_1}\right)\leq \frac{C_1+2}{m}.
    \end{equation*}
    Finally by \cref{e:point-close} we get
    \begin{equation*}
        p_{\mathcal{H}}\left(w_m\frac{(\mathcal{P}-x)_{Y_m'}}{R}+h_m\frac{\pi_{V}^{Y_m}(x_y-x_z)_{V_m}}{R'};\frac{T_{\mtt{i}_{j}\mtt{i}_y\mtt{j}}(K)-T_{\mtt{i}_{j}}(x_z)}{\|A_{\mtt{j}}\|R\|A_{\mtt{i}_{j}}|Y_m\|r_1}\right)\leq \frac{C_1+3}{m}.
    \end{equation*}
    Since this holds for all $y\in\mathcal{B}$, \cref{e:weak-pseudo-tangent} follows.
\end{proof}
The main theorem of this section follows immediately.
\begin{proofref}{t:weakly-dominated-assouad}
    First, recall from \cref{p:projection-assouad-constant} that the function $V\mapsto \dimA \pi_{V^\bot}(K)$ takes constant value $\eta$.
    Therefore, applying \cref{ip:large-tan-slice}, it follows that for all $V\in X_F$, there exists an $E\in\Tan(K)$ and an $x\in\pi_{V^\bot}(E)$ such that
    \begin{equation*}
        \dimH(\pi_{V^\bot}^{-1}(x)\cap E) \geq \dimA K - \eta.
    \end{equation*}
    Conversely, by \cref{t:product-tangent}, for any $V\in X_F$, $E\in\Tan(K)$, and $x\in\pi_{V^\bot}(E)$, there are compact sets $A,B\in\R$ with $\dimlB A=\eta$ and $\dimlB B=\dimA(\pi_{V^{\bot}}^{-1}(x)\cap E)$ such that $A\times B$ is a coarse microset of $K$.
    Therefore by \cref{l:coarse-micro},
    \begin{equation*}
        \begin{aligned}
            \dimA K&\geq \dimA A\times B\geq \dimlB A\times B\geq \dimlB A+\dimlB B\\
            &=\eta+\dimA(\pi_{V^{\bot}}^{-1}(x)\cap E),
        \end{aligned}
    \end{equation*}
    which gives the lower bound.
\end{proofref}

\section{Consequences and examples}\label{s:consequences}
\subsection{Completing the proof of the main result}
In this section, we complete the proof of the main result, \cref{it:tan-slice-form}.

First, let us note the following consequence of \cref{t:weakly-dominated-assouad}, using \cref{ip:large-tan-slice}.
This establishes the main result concerning Assouad dimension and slices of weak tangents.
\begin{corollary}\label{p:sup-attained}
    Let $K$ be a weakly dominated self-affine set.
    Then for any $V\in X_F$ there exists $F\in\Tan(K)$ and $x\in \pi_{V^{\bot}}(F)$ such that
    \begin{equation*}
        \dimH(\pi_{V^{\bot}}^{-1}(x)\cap F)=\dimA(\pi_{V^{\bot}}^{-1}(x)\cap F)=\dimA K-\dimA \pi_{V^{\bot}}(K).
    \end{equation*}
\end{corollary}
\begin{proof}
    Let $V\in X_F$ and denote by $\eta=\sup_{W\in X_F}\dimA\pi_{W^{\bot}}(K)$.
    Recall that by \cref{t:weakly-dominated-assouad},
    \begin{equation*}
        \dimA K=\eta + \sup_{E\in\Tan(K)}\sup_{\substack{V\in X_F\\x\in\pi_{V^{\bot}}(K)}}\dimA(\pi^{-1}(x)\cap F).
    \end{equation*}
    By \cref{ip:large-tan-slice}, there exists $F\in\Tan(K)$ and $x\in\pi(F)$ such that
    \begin{equation*}
        \dimH (\pi_{V^{\bot}}^{-1}(x) \cap F) \geq \dimA K - \dimA \pi_{V^{\bot}}(K).
    \end{equation*}
    Moreover, by \cref{p:projection-assouad-constant}
    \begin{equation*}
        \dimA (\pi_{V^{\bot}}^{-1}(x) \cap F)\leq \sup_{E\in\Tan(K)}\sup_{\substack{V\in X_F\\x\in\pi_{V^{\bot}}(K)}}\dimA(\pi^{-1}(x)\cap F)=\dimA K-\dimA\pi_{V^{\bot}}(K),
    \end{equation*}
    and the claim follows.
\end{proof}
Finally, we use self-affinity and separation conditions to make conclusions about slices of the original set.

We first note the following bound, which is a simple corollary of \cref{t:weakly-dominated-assouad}, by the observation that a weak tangent of a slice of a compact set is contained in a slice of a weak tangent of the set, see e.g.\ \cite[Lemma 4.4]{doi:10.1017/etds.2023.117} for the short formal proof.
This proves the part of \cref{it:tan-slice-form} concerning the upper bound on dimensions of slices with no separation conditions.
\begin{corollary}\label{c:real-slice-ub}
    Let $(T_i)_{i\in\mathcal{I}}$ be a weakly dominated self-affine IFS with attractor $K$.
    Then for all $W\in X_F$ and $x\in\pi_{W^{\bot}}(K)$,
    \begin{equation*}
        \dimA (\pi_{W^{\bot}}^{-1}(x)\cap K)\leq \dimA K-\eta,
    \end{equation*}
    where $\eta$ is the constant value for the map $V\mapsto \dimA \pi_{V^{\bot}}(K)$.
\end{corollary}
Finally, we show that this bound can be upgraded to an equality for some slice under a suitable separation condition which is slightly weaker than the strong separation condition.
\begin{definition}\label{d:wbnc}
    Let $K$ be a self-affine set.
    For $x\in K$ and $r>0$, we let
    \begin{equation*}
        \Phi(x,r)=\{T_{\mtt{i}}\colon \alpha_2(A_{\mtt{i}})\leq r <\alpha_2(A_{\mtt{i}^-})\text{ and }T_{\mtt{i}}(K)\cap B(x,r)\ne\emptyset\}.
    \end{equation*}
    We say that $K$ satisfies the \emph{weak bounded neighbourhood condition (WBNC)} if there exists a constant $M\in\R$ such that
    \begin{equation*}
        \#\Phi(x,r)\leq M,
    \end{equation*}
    for all $x\in K$ and $r>0$.
\end{definition}
This following corollary completes the proof of \cref{it:tan-slice-form}.
\begin{corollary}\label{c:real-slice}
    Suppose $(T_i)_{i\in\mathcal{I}}$ is a weakly dominated IFS satisfying the weak bounded neighbourhood condition with attractor $K$.
    Then there exists $W\in X_F$ and $x\in\pi(K)$ such that
    \begin{equation*}
        \sup_{E\in\Tan(K)}\sup_{\substack{V\in X_F\\x\in\pi_{V^{\bot}}(K)}}\dimA(\pi_{V^{\bot}}^{-1}(x)\cap F) = \dimH(\pi_{W^{\bot}}^{-1}(x)\cap K).
    \end{equation*}
    In particular,
    \begin{equation*}
        \dimA K = \eta + \max_{\substack{V\in X_F\\x\in\pi_{V^{\bot}}(K)}}\dimH(\pi_{V^{\bot}}^{-1}(x)\cap K),
    \end{equation*}
    where $\eta$ is the constant value for the map $V\mapsto \dimA \pi_{V^{\bot}}(K)$.
\end{corollary}
\begin{proof}
    By \cref{p:sup-attained}, choose $F\in\Tan(K)$, $V\in X_F$ and $x\in \pi_{V^{\bot}}(K)$ such that
    \begin{equation*}
        \dimH(\pi_{V^{\bot}}^{-1}(x)\cap F)=\sup_{F\in\Tan(K)}\sup_{\substack{V\in X_F\\x\in\pi_{V^{\bot}}(K)}}\dimA(\pi^{-1}(x)\cap F).
    \end{equation*}
    By \cite[Lemma 3.2]{doi:10.1017/etds.2023.117}, there exists a finite index set $I$ such that
    \begin{equation*}
        F=\bigcup_{i\in I}F_i,
    \end{equation*}
    where each $F_i$ is a compact set, and for each $i\in I$, there exists a linear map $G_i\in \mathfrak{R}(\A^{-1})$, and a point $y_i\in K$ such that $G_i(F_i)+y_i\subset K$.
    Let $j\in I$ be such that $\dimH(\pi_{V^{\bot}}^{-1}(x)\cap F_i)=\dimH(\pi_{V^{\bot}}^{-1}(x)\cap F)$.
    If $G_j$ has rank one, then $\im(G_j)\in X_F$ and $\ker G_j\ne V$ and therefore $G_j(\pi_{V^{\bot}}^{-1}(x)\cap F_i)+y_i\subset \im(G_j)\cap K+y_i$ is bi-Lipschitz equivalent with $\pi_{V^{\bot}}^{-1}(x)\cap F_i$.
    This gives the lower bound
    \begin{equation}\label{e:wbnc-lower-bound}
        \dimH (\pi_{\im(G_j)^{\bot}}^{-1}(y_i)\cap K)\geq \dimH (G_j(\pi_{V^{\bot}}^{-1}(x)\cap F_i)+y_i)\geq \dimH(\pi_{V^{\bot}}^{-1}(x)\cap F).
    \end{equation}

    On the other hand, if $\mathrm{rank}(G_j)=2$, it follows from \cite[Lemma 3.3]{zbl:1452.37039} that $G_jV\in X_F$ and since $G_j$ is globally bi-Lipschitz,
    \begin{equation*}
        \dimH(\pi_{(G_jV)\bot}^{-1}(G_jx+y_i)\cap K)\geq \dimH (G_j(\pi_{V^{\bot}}^{-1}(x)\cap F_i)+y_i)\geq \dimH(\pi_{V^{\bot}}^{-1}(x)\cap F),
    \end{equation*}
    which gives the claim in this case.
\end{proof}

\subsection{Dimensions of tubes and local dimensions}\label{ss:geom-bound}
In \cref{t:weakly-dominated-assouad}, we established a formula for the Assouad dimension of a weakly dominated self-affine set in terms of the Assouad dimension of slices of weak tangents.
We now show that this also enables us to bound the dimensions of tubes.

We introduce some notation: for $V\in\RP^1$ and $z\in \pi_{V}(\R^2)$, we denote the tube with width $2r$ through $z$ in direction $V$ by
\begin{equation*}
    \mathcal{T}_r(V, z) = \pi_{V}^{-1}\bigl(B(z,r)\bigr).
\end{equation*}
We then set
\begin{align*}
    \Delta_{V}(F) = \limsup_{r\to 0}\frac{\sup_{z\in \pi_{V}(K)}\log N_r(\mathcal{T}_r(V,z)\cap F)}{\log(1/r)}.
\end{align*}
Of course, in general,
\begin{equation*}
    \sup_{x\in\pi_{V}(F)}\dimuB (\pi_{V}^{-1}(x)\cap F) \leq \Delta_{V}(F).
\end{equation*}
Conversely, the maximal value of $\Delta_{V}(E)$ over all microsets $E$ is always attained by the Hausdorff dimension of a slice of a weak tangent.
Using \cref{l:h-micro-hausdorff}, the proof is essentially the same as the proof of \cref{ip:large-tan-slice}, so we omit the details.
\begin{lemma}
    Let $F\subset\R^2$ be non-empty and compact and let $V\in\RP^1$ be arbitrary.
    Write $\eta = \sup_{E\in\mathcal{G}_F} \Delta_{V}(E)$.
    Then there exists an $E_0\in\Tan(F)$ and an $x\in\pi_{V}(E_0)$ so that $\mathcal{H}^\eta(\pi_{V}^{-1}(x)\cap E_0) > 0$.
\end{lemma}
Since for all compact sets $F\subset \R^2$ we always have (up to rescaling and translation) $F\in\mathcal{G}_F$ and $\Tan(F)\subset \mathcal{G}_F$, the following corollary of \cref{it:tan-slice-form} is immediate.
\begin{corollary}\label{c:tube-bound}
    Let $(T_i)_{i\in\mathcal{I}}$ be a weakly dominated self-affine IFS with attractor $K$ and let $\eta$ denote the constant value of the map $V\mapsto \dimA \pi_{V^{\bot}}(K)$ for $V\in X_F$.
    Then for all $V\in X_F$,
    \begin{align*}
        \dimA K-\eta &= \max_{E\in\Tan(K)}\Delta_{V^\bot}(E)\geq \Delta_{V^\bot}(K).
    \end{align*}
    If in addition $(T_i)_{i\in\mathcal{I}}$ satisfies the WBNC, then
    \begin{equation*}
        \dimA K - \eta = \max_{V\in X_F}\Delta_{V^\bot}(K) = \max_{V\in X_F}\max_{x\in\pi_{V^\bot}(K)}\dimH(\pi_{V^\bot}^{-1}(x)\cap K).
    \end{equation*}
\end{corollary}
To conclude this section, let us explain how to extend and complete the results of \cite{zbl:1371.28016} using \cref{c:tube-bound}.
Let us begin by recalling the setting and the main result in \cite{zbl:1371.28016}.
Fix numbers $0<\alpha<\beta<1$, and consider the self-affine IFS $(T_i)_{i\in\mathcal{I}}$ given by $T_i(x,y) = (\beta x, \alpha y) + (b_i, a_i)$ where $0\leq b_i\leq 1-\beta$ and $0\leq a_i \leq 1-\alpha$.
We assume that the IFS satisfies the \emph{rectangular open set condition}: $T_i((0,1)^2)\cap T_j((0,1)^2)=\varnothing$ for all $i\neq j$.

Set $m=\#\mathcal{I}$, and let $\mu$ be the unique Borel measure satisfying $\mu(T_{\mtt{i}}((0,1)^2))= m^{-k}$ for $\mtt{i}\in\mathcal{I}^k$.
Equivalently, $\mu$ is just the self-affine measure on the IFS $(T_i)_{i\in\mathcal{I}}$ with equal probabilities.

In our notation, $X_F=\{H\}$ is a singleton and $\pi_{H^\perp}=\pi$ is the orthogonal projection onto the $x$-axis.
Moreover, the rectangular open set condition implies that the IFS satisfies the WBNC.
Let $(S_i(x) = \beta x + b_i)_{i\in\mathcal{I}}$ denote the corresponding projected IFS with attractor $\pi(K)$ and let $\nu=\mu\circ\pi^{-1}$ denote the pushforward measure.
Note that, in particular, $\nu(B(z,r))=\mu(\mathcal{T}_r(H,z))$, for all $z\in\pi(K)$.

We also recall the definition of the \emph{Frostman dimension}:
\begin{equation*}
    \dim_\infty \nu \coloneqq \sup\{t\geq 0:\exists C>0\,\forall r\in(0,1)\,\forall z\in \R \nu(B(x,r)) \leq C r^t\}.
\end{equation*}
Equivalently, $\dim_\infty\nu$ is the slope of the asymptote of the $L^q$ spectrum at $+\infty$, or the infimum of lower local dimensions over all points $z\in\supp\nu$.
The following is a slightly weaker version of the main result of \cite{zbl:1371.28016}.
\begin{theorem}[\cite{zbl:1371.28016}]\label{t:fj-result}
    Let $K$ be the self-affine set corresponding to the IFS $(T_i)_{i\in\mathcal{I}}$ defined above, with parameters $0<\alpha<\beta<1$.
    Let $s$ denote Frostman dimension of $\nu$.
    Then
    \begin{equation*}
        \dimA K \leq \dimA \pi(K) + \frac{\log m\beta^s}{\log(1/\alpha)}.
    \end{equation*}
    Moreover, equality holds if any of the following conditions hold:
    \begin{enumerate}[nl,r]
        \item $\dimB \pi(K) = 1$.
        \item The projected IFS $(S_i)_{i\in\mathcal{I}}$ satisfies the weak separation condition.
        \item The projected IFS $(S_i)_{i\in\mathcal{I}}$ satisfies the exponential separation condition.
    \end{enumerate}
\end{theorem}
Moreover, the conclusion under the assumption of the exponential separation condition relies on the deep work of Shmerkin \cite{zbl:1426.11079}; we also refer the reader to that paper for a precise definition of the exponential separation condition.

We will use \cref{c:tube-bound} to prove that the upper bound for the Assouad dimension in \cref{t:fj-result} is in fact always an equality.
The result will in fact follow from the following simple lemma relating the size of tubes with local dimensions.
Similarly to the other notation, let $\mathcal{T}_r(z)=\mathcal{T}_r(H^\bot, z)$ denote the vertical $r$-tube passing through $x$ and write $\Delta(F) = \Delta_{H^\bot}(F)$.
\begin{lemma}\label{l:local-tube}
    Let $K$ be the self-affine set corresponding to the IFS $(T_i)_{i\in\mathcal{I}}$ defined above, with parameters $0<\alpha<\beta<1$.
    Let $s$ denote the Frostman dimension of $\nu$.
    Then
    \begin{equation*}
        \Delta(K)=\frac{\log m\beta^s}{\log(1/\alpha)}.
    \end{equation*}
\end{lemma}
\begin{proof}
    Let $\varepsilon>0$ be arbitrary.
    Suppose $z\in\pi(K)$ and $r\in(0,1)$.
    Let $n\in\N$ be minimal such that $\alpha^n \leq r$, and let
    \begin{equation}\label{e:lambda}
        \Lambda_n = \{\mtt{i}\in\mathcal{I}^n: S_{\mtt{i}}(\pi(K))\cap B(z,r)\neq\varnothing\}.
    \end{equation}
    Since the IFS satisfies the rectangular open set condition and $\alpha^n\approx r$,
    \begin{equation*}
        N_r(\mathcal{T}_r(z)\cap K) \approx \#\Lambda_n.
    \end{equation*}
    Now fix some rectangle $T_{\mtt{i}}((0,1)^2)$ for $\mtt{i}\in\mathcal{I}^n$ and note that by the self-affinity of $\mu$,
    \begin{equation*}
        \mu\bigl(\mathcal{T}_r(z) \cap T_{\mtt{i}}(K)\bigr) = m^{-n}\mu\bigl(\mathcal{T}_{r \beta^{-n}}(T_{\mtt{i}}^{-1}(z))\bigr) \lesssim m^{-n} (\alpha/\beta)^{n(s-\varepsilon)}.
    \end{equation*}
    Therefore
    \begin{align*}
        \nu(B(z,r)) &\lesssim  N_r(\mathcal{T}_r(z)\cap K)\cdot m^{-n}\cdot (\alpha/\beta)^{n(s-\varepsilon)}\\
                    &\lesssim (\alpha^n)^{-(\Delta(K)+ \varepsilon)}\cdot(\alpha^n)^{\frac{\log m}{\log(1/\alpha)}}(\alpha^n)^{\frac{(s-\varepsilon)\log \beta}{\log(1/\alpha)}} (\alpha^n)^{s-\varepsilon}\\
                    &= (\alpha^n)^{s-\Delta(K) + \frac{\log m\beta^s}{\log(1/\alpha)} - 2\varepsilon - \varepsilon\frac{\log\beta}{\log(1/\alpha)}}.
    \end{align*}
    Since this holds for all $z\in \pi(K)$ and $r\in(0,1)$, by maximality of $s$,
    \begin{equation*}
        \Delta(K) \geq \frac{\log m\beta^s}{\log(1/\alpha)} - 2\varepsilon - \varepsilon\frac{\log\beta}{\log(1/\alpha)}.
    \end{equation*}
    Since $\varepsilon>0$ was arbitrary, the lower bound on $\Delta(K)$ holds.

    Now to obtain the upper bound on $\Delta(K)$, again let $\varepsilon>0$ be arbitrary.
    Fix $z\in\pi(K)$ and $r\in(0,1)$.
    As before, let $n\in\N$ be minimal such that $\alpha^n\leq r$, and let $\Lambda_n$ be defined as in \cref{e:lambda}.
    Then for all $\mtt{i}\in\Lambda_n$,
    \begin{equation*}
        S_{\mtt{i}}(K)\subset B(z, \beta^n + \alpha^n) \subset B(z, 2\beta^n).
    \end{equation*}
    Therefore, recalling that $\#\Lambda_n\approx N_r(\mathcal{T}_r(z)\cap K)$,
    \begin{equation*}
        N_r(\mathcal{T}_r(z)\cap K) \approx \#\Lambda_n \leq m^n\nu\bigl(B(z, 2\beta^n)\bigr) \lesssim m^n \beta^{n(s-\varepsilon)} \approx \left(\frac{1}{r}\right)^{\frac{\log m\beta^{s-\varepsilon}}{\log(1/\alpha)}},
    \end{equation*}
    it follows that
    \begin{equation*}
        \Delta(K) \leq \frac{\log m\beta^{s-\varepsilon}}{\log(1/\alpha)}.
    \end{equation*}
    Since $\varepsilon>0$ was arbitrary, the desired upper bound holds.
\end{proof}
We finally obtain the desired result.
\begin{restatement}{ic:tubes}
    Let $K$ be the self-affine set corresponding to the IFS $(T_i)_{i\in\mathcal{I}}$ defined above, with parameters $0<\alpha<\beta<1$.
    Let $s$ denote the Frostman dimension of $\nu$.
    Then
    \begin{equation*}
        \dimA K = \dimA \pi(K) + \frac{\log m\beta^s}{\log(1/\alpha)}.
    \end{equation*}
\end{restatement}
\begin{proof}
    By \cref{c:tube-bound} (recalling that $X_F$ is a singleton and the IFS satisfies the WBNC) and \cref{l:local-tube},
    \begin{equation*}
        \dimA K - \dimA \pi(K) = \Delta(K) = \frac{\log m\beta^s}{\log(1/\alpha)}
    \end{equation*}
    as claimed.
\end{proof}

\subsection{Dichotomy for the conformal Assouad dimension}\label{ss:conformal}
To conclude our applications, we apply our results to calculate the conformal Assouad dimension of a large class of self-affine sets.
Recall that a function $f\colon X\to Y$ between metric spaces $(X,d)$ and $(Y,\rho)$ is an \emph{$\eta$-quasisymmetry} if
\begin{equation*}
    \frac{\rho(f(x),f(y))}{\rho(f(x),f(z))}\leq \eta\left(\frac{d(x,y)}{d(x,z)}\right),
\end{equation*}
for all $x,y,z\in X$ with $x\ne z$, where $\eta \colon [0,\infty)\to [0,\infty)$ is a homeomorphism.
Quasisymmetries are generalizations of bi-Lipschitz maps which preserve relative sizes of sets, while allowing sets of widely different sizes to be distorted in different ways.
Unlike bi-Lipschitz maps, quasisymmetries can decrease the Assouad dimension which gives rise to the \emph{conformal Assouad dimension} of $X$,
\begin{equation*}
    \CdimA X\coloneqq \inf\left\{\dimA f(X)\colon f\text{ is a quasisymmetry}\right\},
\end{equation*}
which is of course invariant under quasisymmetries.
We say that a set $X$ is \emph{minimal} for conformal Assouad dimension if $\CdimA X=\dimA X$.
We also consider the conformal Hausdorff dimension, which is the same definition except with Hausdorff dimension in place of Assouad dimension.

By combining a strong projection theorem for Assouad dimension due to Orponen with \cref{t:product-tangent}, we are able to prove a dichotomy result for the conformal Assouad dimension for a large class of self-affine sets.
\begin{restatement}{ic:conformal}
    Let $K$ be a weakly dominated and irreducible self-affine set.
    If $\dimA K<1$, then $\CdimA K=0$, and if $\dimA K\geq 1$, then $K$ is minimal for conformal Assouad dimension.
\end{restatement}
\begin{proof}
    If $\dimA K<1$ then \cite[Corollary~5.1.11]{zbl:1201.30002} implies that $\CdimA K=0$.
    Assume therefore that $\dimA K\geq 1$.

    Using \cref{p:sup-attained}, choose $F\in \Tan(K)$, $V\in X_F$ and $x\in\pi_{V^{\bot}}(F)$, such that
    \begin{equation*}
        \dimA K = \dimA \pi_{V^{\bot}}(K)+\dimH(\pi_{V^{\bot}}^{-1}(x)\cap F).
    \end{equation*}
    Note that if $\dimA \pi_{V^{\bot}}(K)=1$ for any $V\in X_F$, since $[0,1]$ is a weak tangent of $\pi_{V^\bot}(K)$, by \cref{p:large-weak-tan} combined with the Lebesgue density theorem, in this case the proof of \cref{t:product-tangent} actually shows that there is a compact set $B$ with $\dimlB B=\dimH(\pi_{V^{\bot}}^{-1}(x)\cap F)$ such that $[0,1]\times B$ is a coarse microset of $K$.
    Moreover, applying \cref{p:large-weak-tan} and passing to a weak tangent again if necessary, we may assume that $\dimH \bigl([0,1]\times B\bigr)=\dimA K$.
    Since quasisymmetries cannot lower the Hausdorff dimension of a product of a compact set with an interval \cite[Proposition~4.1.11]{zbl:1201.30002}, it follows from \cite[Proposition~6.1.5]{zbl:1201.30002}, that
    \begin{equation*}
        \CdimA K\geq \CdimH\bigl([0,1]\times B\bigr)=\dimH\bigl([0,1]\times B\bigr)=\dimA K.
    \end{equation*}

    It therefore remains to show that there exists $V\in X_F$ such that $\dimA \pi_{V^{\bot}}(K)=1$.
    A theorem of Orponen \cite{zbl:1465.28008} shows that for any set $E$ with $\dimA E\geq 1$, we have
    \begin{equation*}
        \dimH\{V\in\RP^1\colon \dimA\pi_{V^{\bot}}(E)<1\}=0.
    \end{equation*}
    It is easy to see (for example using \cite[Lemma 2.7]{arxiv:2107.00983}) that weak domination combined with irreducibility implies strong irreducibility, from which it follows that $\dimH X_F>0$ \cite[Corollary~VI.4.2]{zbl:0572.60001}.
    This gives the claim.
\end{proof}
For \cref{ic:conformal} to hold as stated, irreducibility is necessary: examples of self-affine carpets with Assouad dimension greater than 1 but which are not minimal for conformal Assouad dimension are given in \cite{zbl:1278.37032}.
However, in these examples the conformal Assouad dimension is 0, and we are not aware of any examples of self-affine carpets with Assouad dimension greater than one which are not minimal for the conformal Assouad dimension.

\begin{acknowledgements}
    RA is supported by a doctoral training grant from the Magnus Ehrnrooth Foundation.
    This project began when RA was visiting the University of St Andrews.
    RA thanks the members of the School of Mathematics and Statistics for their hospitality and Vilho, Yrjö and Kalle Väisälä Foundation for funding the visit.

    AR is supported by two grants of Tuomas Orponen, from the European Research Council (ERC) under the European Union's Horizon Europe research and innovation programme (grant agreement No 101087499), and from the Research Council of Finland via the project Approximate incidence geometry, grant no.\ 355453.
    Part of this work was also done while AR was a PhD student at the University of St Andrews, supported by EPSRC Grant EP/V520123/1 and the Natural Sciences and Engineering Research Council of Canada.

    The authors would also like to thank Balázs Bárány, Jonathan Fraser, Antti Käenmäki, Tuomas Orponen, Aleksi Pyörälä and Meng Wu for valuable discussions related to the topics of this paper.
\end{acknowledgements}
\appendix
\section{Proof of the slicing result for diagonal systems}\label{as:diag}
In this section, we include a proof of \cref{t:weakly-dominated-assouad} for weakly dominated diagonal self-affine systems.

Let us set up the required notation.
Fix a planar IFS $(T_i)_{i\in\mathcal{I}}$ where for each $i\in\mathcal{I}$ there are $a_i,b_i\in(0,1)$ and $u_i,v_i\in\R$ so that
\begin{equation*}
    T_i(x,y)= (a_i x+u_i, b_i y+v_i).
\end{equation*}
Let $K$ denote the unique non-empty compact attractor and let $\{S_i(x)=a_i x+u_i\}_{i\in\mathcal{I}}$ denote the projected IFS on the 1\textsuperscript{st} coordinate axis.
Let $\pi\colon\R^2\to\R$ denote the orthogonal projection $\pi(x,y)=x$: then equivalently $S_i$ is the unique map which satisfies $S_i\circ\pi=\pi\circ T_i$.
Of course, $\{S_i\}_{i\in\mathcal{I}}$ is a self-similar IFS with attractor $\pi(K)$.
We refer to such a system as \emph{diagonal}.
We assume that the IFS is \emph{weakly dominated}: in this notation, this means (without loss of generality) that $a_i\geq b_i$ for all $i\in\mathcal{I}$, and $a_i>b_i$ for some $i\in\mathcal{I}$.

We now have the following special case of \cref{t:weakly-dominated-assouad}.
\begin{theorem}\label{t:diag-proof}
    Let $(T_i)_{i\in\mathcal{I}}$ be a diagonal weakly dominated IFS with attractor $K$.
    Then
    \begin{equation*}
        \dimA K = \dimA\pi(K)+\sup_{E\in\Tan(K)}\sup_{x\in\pi(E)}\dimA(\pi^{-1}(x)\cap E).
    \end{equation*}
\end{theorem}
\begin{proof}
    We recall by \cref{p:large-weak-tan} that it suffices to prove the lower bound on $\dimA K$.

    Let $E\in\Tan(K)$ and $x\in\pi(E)$ be arbitrary.
    For notational simplicity, write $\eta=\dimA\pi(K)$ and $\beta=\dimA(\pi^{-1}(x)\cap E)$.
    We will show that $\dimA K\geq\eta+\beta$ by constructing a coarse microset $G$ with $\dimlB G\geq\eta+\beta$.

    First, applying \cref{c:weak-tan}, get a sequence of dyadic cubes $(P_m)_{m=1}^\infty$ with $P_m\in\mathcal{D}_{k_m}$ satisfying $\lim_{m\to\infty}k_m=\infty$ such that
    \begin{equation}\label{e:proj-cover}
        N_{2^{-n}}(\psi_{P_m}(\pi(K))\cap[0,1])\gtrsim 2^{n\left(\eta-\frac{1}{m}\right)}
    \end{equation}
    for all $n\in\Z$ with $0\leq n\leq m$.

    Now, fix an $m\in\N$.
    We begin by pigeonholing a good set of cylinders.
    Let $M\in\N$ be sufficiently large so that
    \begin{equation*}
        M\geq  \bigl(m 2^{k_m}\bigr)^2\cdot (1-a_{\min})\cdot \bigl(8 a_{\min}^{-1}\cdot\diam \pi(K)\bigr).
    \end{equation*}
    First, by definition of the Assouad dimension, get $y\in \pi^{-1}(x)\cap E$, $0<r_0\leq R_0<1$ with $R_0\geq 2^{m(m+k_m)}\cdot r_0$, and a finite set of points $\mathcal{A}_0\subset E_x$ which are $6r_0$-separated with
    \begin{equation*}
        \#\mathcal{A}_0\geq 4M\left(\frac{R_0}{r_0}\right)^{\beta-\frac{1}{m}}.
    \end{equation*}
    Since $E\in\Tan(K)$, there exist $z\in K$ and $\lambda \geq 1$ so that
    \begin{equation*}
        d_{\mathcal{H}}\bigl(\lambda (K-z)\cap B(0,1), E\bigr)\leq r_0.
    \end{equation*}
    Let $\ell_{m,1}\in\N\cup\{0\}$ be maximal and $\ell_{m,2}\in\N\cup\{0\}$ be minimal so that
    \begin{equation*}
        2^{-\ell_{m,2}} \leq \lambda^{-1} r_0 \leq  \lambda^{-1} R_0 \leq 2^{-\ell_{m,1}}.
    \end{equation*}
    Write $R = 2^{-\ell_{m,1}}$ and $r= 2^{-\ell_{m,2}}$.
    Then for each $a\in\mathcal{A}_0$, there is a $\xi(a)\in E$ such that $d(\lambda^{-1} a + z,\xi(a))\leq r_0 \lambda^{-1}$.

    Let $\mathcal{A}=\{\xi(a):a\in\mathcal{A}_0\}$.
    We observe four key properties of $\mathcal{A}$.
    \begin{enumerate}[nl]
        \item Since the points in $\mathcal{A}_0$ are $6r_0$-separated, the points in $\mathcal{A}$ are $2r$-separated.
        \item By the condition on $R_0/r_0$, we have $\ell_{m,2}-\ell_{m,1}\geq m(m+k_m)$.
        \item Since $\beta-1/m \leq 1$,
            \begin{equation*}
                \#\mathcal{A} = \#\mathcal{A}_0 \geq 4M\left(\frac{R_0}{r_0}\right)^{\beta-\frac{1}{m}} \geq M \left(\frac{R}{r}\right)^{\beta - \frac{1}{m}}.
            \end{equation*}
        \item Since $\mathcal{A}_0$ is a subset of a vertical slice, $\mathcal{A}$ is a subset of a vertical tube of width $2\lambda^{-1} r_0\leq 4r$.
    \end{enumerate}

    Now for each $y\in \mathcal{A}$, let $\mtt{i}_y\in\mathcal{I}^*$ be such that $a_{\mtt{i}}\leq r<a_{\mtt{i}^-}$ and $(x,y_j)\in T_{\mtt{i}}(K)$.
    Since the points in $\mathcal{A}$ are $2r$-separated, the images $T_{\mtt{i}_y}(K)$ are disjoint for distinct $y\in \mathcal{A}$.
    Moreover, writing $S_{\mtt{i}_y}(x)= u(y) x+v(y)$, for all $y,z\in \mathcal{A}$, since the points in $\mathcal{A}$ lie in a vertical tube of width $4r$,
    \begin{equation*}
        \frac{u(y)}{r}\in[a_{\min},1]\quad\text{and}\quad\left\lvert\frac{v(y)-v(z)}{r}\right\rvert\leq 4\cdot a_{\min}^{-1}\cdot\diam\pi(K).
    \end{equation*}
    Thus by the pigeonhole principle and the choice of $M$, get $u_m\in[a_{\min},1]$, $v_m\in[0,a_{\min}^{-1}\cdot\diam\pi(K)]$ and $\mathcal{B}\subset \mathcal{A}$ where $\#\mathcal{B}\geq \#\mathcal{A}/M$ such that for all $y\in \mathcal{B}$,
    \begin{equation}\label{e:cyl-pigeonhole}
        \left\lvert\frac{u(y)}{r}-u_m\right\rvert\leq \frac{1}{m 2^{k_m}}\quad\text{and}\quad\left\lvert\frac{v(y)}{r}-v_m\right\rvert\leq \frac{1}{m 2^{k_m}}.
    \end{equation}
    Since the points in $\mathcal{B}$ are $2r$-separated, each $y\in \mathcal{B}$ intersects a distinct vertical dyadic interval of width $2^{-\ell_{m,2}}$.
    On the other hand, $\mathcal{B}$ intersects at most $2$ dyadic intervals of width $2^{-\ell_{m,1}}$.
    Thus pigeonholing again, get $\mathcal{E}_m\subset \mathcal{B}$ so that
    \begin{equation*}
        \#\mathcal{E}_m\geq \left(\frac{R}{r}\right)^{\beta-\frac{1}{m}}
    \end{equation*}
    and moreover $\mathcal{E}_m$ is contained in a single dyadic interval $Q_0$ of width $2^{-\ell_{m,1}}$.
    Finally, since $\ell_{m,2}-\ell_{m,1}\geq m(m+k_m)$, by applying \cref{l:h-micro} to the set $\psi_{Q_0}(\mathcal{E}_m)\cap[0,1]$ with $s=\beta-2/m$, $t=\beta-1/m$, $\ell=m$, and $k=m+k_m$, get a dyadic interval $Q_m\in\mathcal{D}_{\ell_m}$ with $\ell_{m,1}\leq \ell_m\leq \ell_{m,2}-m-k_m$ such that
    \begin{equation}\label{e:slice-cover}
        N_{2^{-n}}\left(\psi_{Q_m}(\mathcal{E}_m)\cap[0,1]\right)\gtrsim 2^{n\left(\beta-\frac{2}{m}\right)}
    \end{equation}
    for all $n\in\Z$ with $0\leq n\leq m$.

    We now construct a coarse microset $G$ with $\dimlB G\geq \eta+\beta$.
    First, for each $m\in\N$, set
    \begin{equation*}
        g_m(x) = r(u_m x + v_m).
    \end{equation*}
    In light of \cref{e:cyl-pigeonhole}, one should think of the function $g_m$ as an approximation for the functions $S_{\mtt{i}_y}$ for $y\in\mathcal{E}_m$.
    Next, we set
    \begin{equation*}
        F_m= g_m(P_m\cap \pi(K)) \times (Q_m\cap\mathcal{E}_m)\quad\text{and}\quad h_m=2^{-k_m-\ell_{m,2}}.
    \end{equation*}
    Note that $g_m(P_m)$ is an interval of width $u_m\cdot h_m$ and $Q_m$ is an interval of width $2^{-\ell_m}\geq 2^m h_m$.
    Since $(T_i)_{i\in\mathcal{I}}$ is weakly dominated, get $i_0\in\mathcal{I}$ such that $\kappa\coloneqq b_{i_0}/a_{i_0}\in(0,1)$.
    For each $m\in\N$, let $j_m\in\Z$ be maximal such that
    \begin{equation*}
        \kappa^{j_m} 2^{-\ell_m}\geq h_m.
    \end{equation*}
    Observe that $\lim_{m\to\infty}j_m=\infty$.

    Now passing to a subsequence, we may assume that $j_m\geq 0$ for all $m\in\N$ and that the limits $u\coloneqq\lim_{m\to\infty}u_m$ and $v\coloneqq\lim_{m\to\infty}v_m$ exist.
    Set
    \begin{equation*}
        E_m= T_{i_0}^{j_m}(F_m)\qquad\text{and}\qquad w_m=a_{i_0}^{j_m}h_m.
    \end{equation*}
    Of course, $E_m$ is contained in the rectangle
    \begin{equation*}
        \Delta_m= T_{i_0}^{j_m}\left(g_m(P_m)\times Q_m\right)
    \end{equation*}
    which has width $u_m\cdot w_m$ and height $\kappa^{-1}w_m$.
    Let $f_m$ denote the orientation-preserving diagonal affine map satisfying $f_m(\Delta_m)=[0,1]^2$.
    Passing again to a subsequence if necessary, we may assume that the limit
    \begin{equation}\label{e:micro-def}
        G=\lim_{m\to\infty}g_m(E_m)
    \end{equation}
    exists.

    We first show that $G$ is a coarse microset of $K$.
    Let $m\in\N$ and fix $y\in Q_m\cap\mathcal{E}_m$.
    Consider the horizontal strips
    \begin{align*}
        X_{m,y}&\coloneqq T_{i_0}^{j_m}\left(g_m(P_m\cap \pi(K)) \times \{y\}\right),\\
        Y_{m,y}&\coloneqq S_{i_0}^{j_m}\circ S_{\mtt{i}_y}(P_m\cap \pi(K))\times\{T_{i_0}^{j_m}(y)\}.
    \end{align*}
    By \cref{e:cyl-pigeonhole},
    \begin{equation*}
        d_{\mathcal{H}}\left(X_{m,y}, Y_{m,y}\right)\lesssim \frac{1}{m} a_{i_0}^{j_m}2^{-k_m-\ell_{m,2}}=\frac{1}{m} w_m.
    \end{equation*}
    Moreover, since $y\in T_{\mtt{i}_y}(K)$ and the cylinder $T_{\mtt{i}_y}(K)$ has height $\lesssim a_{\mtt{i}_y}$,
    \begin{equation*}
        d_{\mathcal{H}}\bigl(Y_{m,y}, T_{i_0}^{j_m}\circ T_{\mtt{i}_y}(K)\bigr)\lesssim b_{i_0}^{j_m}2^{-\ell_{m,2}}\lesssim \kappa^{j_m} w_m.
    \end{equation*}
    Since this holds for all $y\in Q_m\cap\mathcal{E}_m$,
    \begin{equation*}
        p_{\mathcal{H}}\left(E_m; T_{i_0}^{j_m}(K)\right)\lesssim\left(\frac{1}{m}+\kappa^{j_m}\right)\cdot w_m.
    \end{equation*}
    Since the rectangle $\Delta_m$ has height and width approximately $w_m$, it follows that $G$ is a coarse microset of $K$.

    We now conclude the proof of the lower bound.
    By \cref{e:proj-cover} and \cref{e:slice-cover} and the definition of $g_m$, for all $1>r\geq 2^{-m}$,
    \begin{equation*}
        N_r(g_m(E))\gtrsim \left(\frac{1}{r}\right)^{\eta+\beta-\frac{3}{m}}.
    \end{equation*}
    Thus by \cref{e:micro-def} and \cref{l:lower-box}, $\dimlB G\geq \eta+\beta$, so by \cref{l:coarse-micro}, $\dimA K\geq \eta+\beta$.
    Since $E\in\Tan(K)$ and $x\in\pi(E)$ were arbitrary, the result follows.
\end{proof}
\end{document}